\documentclass[a4paper,12pt]{article}
\usepackage{amsmath,amssymb,amsthm,epic,eepic,eepicemu,graphicx}
\usepackage[english]{babel}

\title{Geometry of Yang--Baxter maps:\\
       pencils of conics and quadrirational mappings}
\author{V.E.\,Adler$^{*,a}$ \and 
        A.I.\,Bobenko$^{*,b}$ \and 
	Yu.B.\,Suris$^{*,s}$}

\date{\today} 

\def\a{\alpha}
\def\b{\beta}
\def\g{\gamma}
\def\d{\delta}

\def\ep{\epsilon}
\def\l{\lambda}
\def\k{\varkappa}

\newcommand{\cC}{{\cal C}}
\newcommand{\cD}{{\cal D}}
\newcommand{\cF}{{\cal F}}
\newcommand{\cG}{{\cal G}}
\newcommand{\cM}{{\cal M}}
\newcommand{\cT}{{\cal T}}
\newcommand{\cX}{{\cal X}}
\newcommand{\Rho}{{\rm P}}

\def\F{\bar F}      
\def\R{\bar R}

\def\I{{\rm I}}
\def\II{{\rm II}}
\def\III{{\rm III}}
\def\IV{{\rm IV}}
\def\V{{\rm V}}

\def\Complex{\mathbb C}
\def\Real{\mathbb R}
\def\CP{\mathbb C\mathbb P}

\def\hat{\widehat}
\def\tot{\leftrightarrow}

\def\({\left(}
\def\){\right)}
\def\<{\langle}
\def\>{\rangle}

\def\To#1{\overset{#1}\longmapsto}

\def\const{\mathop{\rm const}}
\def\itbf{\itshape\bfseries}

\newtheorem{theorem}{Theorem}
\newtheorem{lemma}[theorem]{Lemma}
\newtheorem{proposition}[theorem]{Proposition}

\newtheorem{definition}[theorem]{Definition}

\newtheorem{Example}{Example}[section]

\begin{document}
\def\figurename{Fig.}
\let\oldf=\thefootnote \renewcommand{\thefootnote}{}

\footnotetext{$^*$ 
 Institut f\"ur Mathematik, Technische Universit\"at Berlin,
 Str.~des 17.~Juni 136, 10623 Berlin, Germany.}

\footnotetext{$^{a)}$ On leave from: 
 Landau Institute for Theoretical Physics, 
 Institutsky pr. 12, 142432 Chernogolovka, Russia. 
 E-mail: {\tt adler@itp.ac.ru}.
 Supported by the Alexander von Humboldt Stiftung and by the RFBR grant 
 02-01-00144.}

\footnotetext{$^{b)}$ E-mail: {\tt bobenko@math.tu-berlin.de}.
 Partially supported by the SFB 288 ``Differential Geometry and
 Quantum Physics'' and the DFG Research Center ``Mathematics for Key
 Technologies'' (FZT 86) in Berlin.}

\footnotetext{$^{s)}$ E-mail: {\tt suris@sfb288.math.tu-berlin.de}.
 Supported by the SFB 288 ``Differential Geometry and
 Quantum Physics''.}

\let\thefootnote=\oldf

\maketitle

\begin{abstract}
Birational Yang-Baxter maps (`set-theoretical solutions of the Yang-Baxter 
equation') are considered. A birational map $(x,y)\mapsto(u,v)$ 
is called quadrirational, if its graph is also a graph of a birational map
$(x,v)\mapsto(u,y)$.
We obtain a classification of quadrirational maps on $\CP^1\times\CP^1$,
and show that all of them satisfy the Yang-Baxter equation. 
These maps possess a nice geometric interpretation in terms of linear
pencil of conics, the Yang-Baxter property being interpreted as a new
incidence theorem of the projective geometry of conics. 
\end{abstract}

\paragraph{Keywords:} Yang-Baxter map, Yang-Baxter equation, 
set-theoretical solution, 3D-consistency, quadrirational map

\eject

\section{Introduction}\label{s:intro}

Our major object of interest in this paper are {\itbf Yang--Baxter maps}, 
i.e. invertible maps satisfying the Yang--Baxter equation
\begin{equation}\label{YBE}
 R_{23}\circ R_{13}\circ R_{12}=R_{12}\circ R_{13}\circ R_{23}\,.
\end{equation}
Here $R_{ij}:\cX^3\to\cX^3$ acts on the $i$-th and $j$-th factors of the
cartesian product $\cX^3$ and is identical on the rest one. The maps $R_{ij}$,
considered as mappings $\cX\times\cX\to\cX\times\cX$, may not coincide with
each other, for example they can correspond to different values of some
parameters.

This notion was introduced in \cite{D} under the name of {\it set-theoretical 
solutions of the Yang--Baxter equation}; a number of examples for the case 
when $\cX$ is a finite set was found in \cite{Hi}, and classification results 
for this case were achieved in \cite{ESS}, see also \cite{LYZ}. It 
turned out that the following non-degeneracy notion is natural and useful 
for such maps. Let 
\begin{equation}
R:(x,y)\mapsto (u,v)=(u(x,y),v(x,y)); 
\end{equation}
this map is called {\it non-degenerate} if for any fixed $x\in\cX$, the map 
$v(x,\cdot):\cX\mapsto\cX$ is bijective, and for any fixed $y\in\cX$, the map 
$u(\cdot,y):\cX\mapsto\cX$ is bijective. 

Important examples of Yang--Baxter maps were found in 
\cite{A93,NY,KNY,E,V02,GV}. It was pointed out in \cite{E} that the 
natural ambient category for these examples is not that of sets, but that of
irreducible algebraic varieties, so that the role of morphisms is played
by birational isomorphisms rather than by bijections. So, in this setting
it is supposed that $R$ is a birational isomorphism of $\cX\times\cX$ to
itself. Then the rational maps $v(x,\cdot):\cX\mapsto\cX$ and 
$u(\cdot,y):\cX\mapsto\cX$ are well defined for generic $x$, resp. $y$,
and $R$ is called non-degenerate, if both these maps are birational 
isomorphisms of $\cX$ to itself.  

In \cite{V02} the term ``Yang--Baxter maps'' was proposed instead of
``set-theoret\-ical solutions'', and various notions of {\it integrability}
were studied. In particular, commuting monodromy maps were constructed and Lax
representations were discussed. A general construction of Lax representations 
was given subsequently in \cite{SV}. In the present paper we want to support 
the change of terminology proposed in \cite{V02}, since we consider the term
``Yang--Baxter maps'' much more suitable and adequate. Correspondingly, we 
will say {\it birational Yang-Baxter maps} instead of ``rational 
set-theoretical $R$-matrix'' (which was proposed initially in \cite{E}).
Moreover, since the term ``non-degenerate'' is heavily overloaded in the
mathematical literature, we propose a new term for non-degenerate (in the above
sense) birational Yang-Baxter maps, namely we will call them {\itbf
quadrirational}. We find this notation highly suggestive, in particular because
of the following way to visualize the Yang-Baxter equation (\ref{YBE}).  Assign
elements of $\cX$ (``fields'') to the edges of a quadrilateral, so that  the
arrow from the down-left to the up-right corner of the left square on 
Fig.~\ref{fig:maps} encodes the map $R:(x,y)\mapsto(u,v)$. Then the opposite
arrow naturally encodes the inverse map $R^{-1}:(u,v)\mapsto(x,y)$, and the
quadrirationality (non-degeneracy in the previous terminology) means the
existence of the rational maps $\R:(u,y)\mapsto(x,v)$ and 
$\R^{-1}:(x,v)\mapsto(u,y)$, which are encoded by two arrows on the right
square on Fig.~\ref{fig:maps}. We call $\R,\R^{-1}$ the {\it companion maps}
for $R$. 

\def\tmpg{\path(0,0)(0,100)(100,100)(100,0)(0,0)
 \put(47,-13){$x$}\put(-13,47){$y$}\put(47,105){$u$}\put(105,47){$v$}}
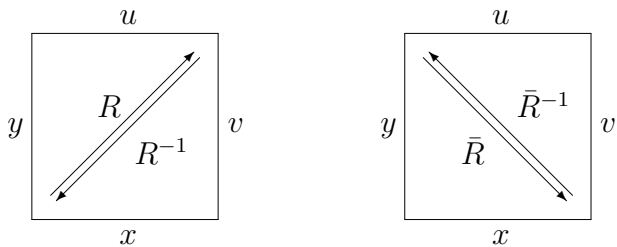
\begin{figure}[htbp]
\begin{center}
\setlength{\unitlength}{0.06em}
\begin{picture}(300,140)(0,-15)
  \tmpg \put(10,13){\vector(1,1){77}}   \put(35,55){$R$}
        \put(90,87){\vector(-1,-1){77}} \put(55,30){$R^{-1}$} 
 \put(200,0){
  \tmpg \put(90,13){\vector(-1,1){77}}  \put(30,30){$\R$} 
        \put(10,87){\vector(1,-1){77}}  \put(60,55){$\R^{-1}$}}
\end{picture}
\caption{A map on $\cX\times\cX$, its inverse and its companions}
\label{fig:maps}
\end{center}
\end{figure}

The Yang--Baxter equation (\ref{YBE}) can be depicted as consistency of
the maps $R$ attached to 6 facets of a cube. The right-hand side of (\ref{YBE})
corresponds to the chain of maps along three front faces of the cube on
Fig.~\ref{fig:YB}:
\[
 (y,z)     \To{R_{23}} (y_3,z_2),    \quad 
 (x,z_2)   \To{R_{13}} (x_3,z_{12}), \quad
 (x_3,y_3) \To{R_{12}} (x_{23},y_{13})
\]
while its left-hand side corresponds to the analogous chain along the rear
faces (it is supposed that opposite faces of the cube carry the same maps). 
Eq.~(\ref{YBE}) assures that two ways of obtaining
$(x_{23},y_{13},z_{12})$ from the initial data $(x,y,z)$ lead to the same
results.

\def\ini(#1,#2){\put(#1,#2){\circle*{5}}}
\def\tmpg{ 
 \begin{minipage}[t]{150pt}
  \begin{picture}(200,150)(-80,-50)   
   \ini(-20,-30) 
   \dashline{5}(-60,-20)(0,0)(0,100)\dashline{5}(0,0)(80,-20)
   \path(-60,80)(-60,-20)(20,-40)(80,-20)(80,80)(0,100)
        (-60,80)(20,60)(80,80) \path(20,60)(20,-40) }
\begin{figure}[htbp]\setlength{\unitlength}{0.07em}
\begin{center}
 \tmpg \ini(50,-30) \ini(80,30) 
   \put(75,-12){\vector(-3,4){50}} 
   \put(13,-28){\vector(-2,3){65}} 
   \put(17, 65){\vector(-1,2){14}}
   \put(-25,-43){$x$}   \put( 36,-20){$x_2$}
   \put(-30, 60){$x_3$} \put( 35, 98){$x_{23}$}
   \put( 45,-43){$y$}   \put(-28,-18){$y_1$} 
   \put( 45, 60){$y_3$} \put(-40, 98){$y_{13}$}
   \put(85,27){$z$}     \put(23,7){$z_2$}  
   \put(-15,40){$z_1$}  \put(-80,25){$z_{12}$}
   \put(-38,10){$R_{13}$}\put(15,75){$R_{12}$}\put(48,30){$R_{23}$}
 \end{picture}
  \caption{Right-hand side of Yang-Baxter equation}\label{fig:YB}
 \end{minipage}\qquad
 \tmpg \ini(-30,-10) \ini(-60,30) 
   \dashline{5}(-55,-10)(-10,80)\put(-8,84){\vector(1,2){4}} 
   \dashline{5}(-45,-20)(60,-20)\put(65,-20){\vector(1,0){4}} 
   \put(-54,-15){\vector(1,1){70}}
   \put(-25,-43){$x$}   \put( 40,-5){$x_2$}
   \put(-35, 63){$x_3$} \put( 35, 96){$x_{23}$}
   \put(-23,-10){$y$}   \put( 45,-43){$y_1$}  
   \put(-40, 96){$y_3$} \put( 45, 60){$y_{13}$} 
   \put(-73,25){$z$}    \put(2,73){$z_2$}
   \put( 23,15){$z_1$}  \put(83,27){$z_{12}$}     
   \put(-25,6){$F_{13}$}\put(-52,40){$F_{23}$}\put(-8,-27){$F_{12}$}
 \end{picture}
 \caption{3D consistency}\label{fig:3D}
 \end{minipage}
\end{center} 
\end{figure}

An alternative way to express this idea is known under the name of {\em 3D
consistency} property (cf. \cite{ABS}). It uses another setting of initial 
data on the edges, as shown on Fig.~\ref{fig:3D}. Also, to distinguish 
between these two notions, we use a different notation for the maps ($F_{ij}$ 
instead of $R_{ij}$) on this figure. Applying the maps to initial data one 
obtains 6 new values at the first step:
\[
 (x,y)\To{F_{12}}(x_2,y_1),\quad 
 (x,z)\To{F_{13}}(x_3,z_1),\quad
 (y,z)\To{F_{23}}(y_3,z_2),
\]
and applying the maps once more, one gets two values for each of the fields
$x_{23},y_{13},z_{12}$. We say that the maps $F_{ij}$ are 3D consistent if 
the two values for each of these fields coincide (for arbitrary initial data).

Comparing the two pictures, we see that (at least in the situation when the
maps are quadrirational) the two notions are equivalent, if $F_{ij}=\R_{ij}$. 

As pointed out above, all known birational Yang-Baxter maps turn out to possess
a more strong property of being quadrirational. In our opinion, the notion of a
quadrirational map is of a great interest for its own, and our results indicate
that it is intimately related with the 3D consistency (although the precise
relation still has to be uncovered). In any case, a more detailed understanding
of the nature of quadrirational maps seems  to be a necessary pre-requisite for
any attempt to classify 3D consistent maps (or, equivalently, Yang-Baxter
maps). In the present paper we furnish this  task in the case $\cX=\CP^1$. 

We start in Sect.~\ref{s:geo} with a precise definition of quadrirationality, 
and then give a nice construction, related to the geometry of conics, which
leads to a large class of quadrirational maps. It allows also to prove the 3D
consistency of the maps corresponding to a triple of conics belonging to one
linear pencil. (See Theorem \ref{th:greeks}; this typical incidence statement
of projective geometry is probably known, however we have not succeeded in
finding predecessors.) In other words, linear pencils of conics carry a
remarkable family of novel Yang--Baxter maps. Upon a rational parametrization
of all encountered conics, these Yang--Baxter maps act on $\CP^1\times\CP^1$.

Then, we proceed with the task of classifying quadrirational maps on
$\CP^1\times\CP^1$. An explicit description of all such maps  is given in
Sect.~\ref{s:BD}. They are {\itbf bi-M\"obius}, i.e. have the form
\begin{equation}\label{map}
  F:\;\; u=\frac{a(y)x+b(y)}{c(y)x+d(y)}\,,\quad
         v=\frac{A(x)y+B(x)}{C(x)y+D(x)}\,,
\end{equation}
where $a(y),\ldots,d(y)$ are polynomials in $y$, while $A(x),\ldots,D(x)$ are
polynomials in $x$. It turns out that there exist three subclasses of such
maps, depending on the degree of the polynomials involved; they are denoted by
pair of numbers as [1:1], [1:2] and [2:2], corresponding to the highest
degrees of the coefficients of both fractions in (\ref{map}). The most rich and
interesting subclass is [2:2]; it is analyzed in detail, while the two simpler
ones are briefly discussed in Appendix \ref{a:12}. 

By classifying quadrirational maps on $\CP^1\times\CP^1$, it is natural  to
factor out the action of the group $({\cal M}\ddot{o}b)^4$ of M\"obius 
transformations acting independently on each field $x,y,u,v$ (recall that 
M\"obius transformations are birational isomorphisms of $\CP^1$). The
construction of Sect.~\ref{s:BD} is not well suited for this purpose. This is
done in Sect.~\ref{s:sing}--\ref{s:class} with the help of different techniques
based on the analysis of the singularities of the maps, which seems to be an
adequate language and admits further generalizations. We obtain the following
result.

\begin{theorem}\label{Th: intro 1}
Any quadrirational map on $\CP^1\times\CP^1$ of subclass [2:2] is equivalent,
under some change of variables from $({\cal M}\ddot{o}b)^4$, to exactly one
of the following five maps: 
\begin{align}
\label{F1}\tag{$F_\I$}
  u&= \a yP,& 
  v&= \b xP,&
  P&= \frac{(1-\b)x+\b-\a+(\a-1)y}
           {\b(1-\a)x+(\a-\b)yx+\a(\b-1)y},& \\
\label{F2}\tag{$F_\II$}
  u&= \frac{y}{\a}\,P,& 
  v&= \frac{x}{\b}\,P,& 
  P&= \frac{\a x-\b y+\b-\a}{x-y},& \\
\label{F3}\tag{$F_\III$}
  u&= \frac{y}{\a}\,P,& 
  v&= \frac{x}{\b}\,P,& 
  P&= \frac{\a x-\b y}{x-y},& \\
\label{F4}\tag{$F_\IV$}
  u&= yP,& 
  v&= xP,& 
  P&= 1+\frac{\b-\a}{x-y},& \\ 
\label{F5}\tag{$F_\V$}
  u&= y+P,& 
  v&= x+P,& 
  P&= \frac{\a-\b}{x-y},&
\end{align}
with some suitable constants $\a,\b$.
\end{theorem}

Remarkably, these five maps describe five particular cases of the geometrical
construction of Sect.~\ref{s:geo}. To mention several further features of these
maps: each one of them is an involution and coincides with its companion maps,
so that all four arrows on Fig.~\ref{fig:maps} are described by the same
formulae. Moreover, these maps have the following additional symmetry:
\[
  x\tot y,\; u\tot v,\; \a\tot\b .
\]
Note also that these maps come with the intrisically built-in parameters
$\a,\b$. Neither their existence nor a concrete dependence on parameters are
pre-supposed in Theorem \ref{Th: intro 1}. The geometric interpretation of 
these parameters will be given in Sect.~\ref{s:class}. It will be shown that 
in the formulas above the parameter $\a$ is naturally assigned to the edges 
$x,u$, while $\b$ is naturally assigned to the edges $y,v$.

The most remarkable fact about the maps \ref{F1}\,--\ref{F5} is their 3D
consistency. For any $\cT=\I,\II,\III,\IV$ or $\V$, denote the corresponding
map $F_\cT$ of Theorem \ref{Th:  intro 1} by  $F_{\cT}(\a,\b)$, indicating the
parameters explicitly. 

\begin{theorem}\label{Th: intro 2}
For any $\a_1,\a_2,\a_3\in\Complex$, the maps $R_{ij}=F_{\cT}(\a_i,\a_j)$ 
(acting nontrivially on the $i$-th and the $j$-th factors of $(\CP^1)^3$)
satisfy the Yang--Baxter equation (\ref{YBE}). Analogously, the maps
$F_{ij}=F_{\cT}(\a_i,\a_j)$ are 3D consistent.
\end{theorem}

In other words, quadrirational maps on $\CP^1\times\CP^1$ of the subclass
[2:2] almost automatically (that means, in a suitable normalization)
are Yang--Baxter maps. Similar statements hold also for other two subclasses
[1:1] and [1:2], so that in the scalar case the properties  of being
quadrirational and of being 3D consistent are related very closely. It
would be intriguing to establish such a relation for more general algebraic
varieties, e.g. for $\CP^n\times\CP^n$.

\section{Geometry of conics and quadrirational \\ maps}\label{s:geo}

We start with giving a slightly more general definition of quadrirational maps.
Let $\cX_i$ $(i=1,2)$ be irreducible algebraic varieties over $\Complex$. 
Consider a rational map $F:\cX_1\times\cX_2\to\cX_1\times\cX_2$. It will be 
identified with its graph -- an algebraic variety
$\Gamma_F\subset\cX_1\times\cX_2\times\cX_1\times\cX_2$.

\begin{definition}
A map $F$ is called {\itbf quadrirational}, if for any fixed pair
$(X,Y)\in\cX_1\times\cX_2$ except maybe some closed subvarieties of 
codimension $\ge 1$, the variety $\Gamma_F$ intersects each one of the sets
$\{X\}\times\{Y\}\times\cX_1\times\cX_2$, 
$\cX_1\times\cX_2\times\{X\}\times\{Y\}$,
$\cX_1\times\{Y\}\times\{X\}\times\cX_2$ and
$\{X\}\times\cX_2\times\cX_1\times\{Y\}$ exactly once, i.e. if
$\Gamma_F$ is a graph of four rational maps 
$F,F^{-1},\F,\F^{-1}:\cX_1\times\cX_2\mapsto\cX_1\times\cX_2$.
\end{definition}

This definition is immediately applicable to the following map. Consider a pair
of nondegenerate conics $Q_1$, $Q_2$ on the plane $\CP^2$, so that both $Q_i$
are irreducible algebraic varieties isomorphic to $\CP^1$.  Take $X\in Q_1$,
$Y\in Q_2$, and let $\ell=\overline{XY}$ be the line through $X,Y$
(well-defined if $X\neq Y$). Generically, the line $\ell$ intersects  $Q_1$ at
one further point $U\neq X$, and intersects $Q_2$ at one further  point $V\neq
Y$. This defines a map $\cF:(X,Y)\mapsto(U,V)$, see Fig.~\ref{fig:greeks} for
the $\Real^2$ picture.

\begin{figure}[htbp]
\begin{center}\includegraphics[width=8cm]{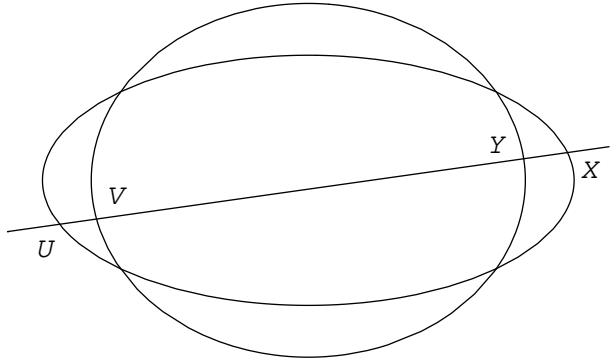}\end{center}
\caption{A quadrirational map on a pair of conics}\label{fig:greeks}
\end{figure}

\begin{proposition}
The map $\cF:Q_1\times Q_2\mapsto Q_1\times Q_2$ is quadrirational. It is an
involution and moreover coincides with its both companions. Intersection 
points $X\in Q_1\cap Q_2$ correspond to the singular points $(X,X)$ of the 
map $\cF$.
\end{proposition}
\begin{proof} 
Intersection of a conic with a line is described by a quadratic equation. 
Knowing one root of a quadratic equation allows us to find the second
one as a rational function of the input data. 
\end{proof}

Generically, two conics intersect in four points, however, degeneracies
can happen. Five possible types $\I-\V$ of intersection of two conics are 
described in detail in \cite{Berger}:
\begin{itemize}\setlength{\itemsep}{0pt}
\item[\I:] four simple intersection points;
\item[\II:] two simple intersection points and one point of tangency;
\item[\III:] two points of tangency;
\item[\IV:] one simple intersection point and one point of the second order 
tangency;
\item[\V:] one point of the third order tangency.
\end{itemize}
Using rational parametrizations of the conics:
\[
\CP^1\ni x\mapsto X(x)\in Q_1\subset\CP^2,\quad {\rm resp.}\quad
\CP^1\ni y\mapsto Y(y)\in Q_2\subset\CP^2,
\]
it is easy to see that $\cF$ pulls back to the map
$F:(x,y)\mapsto(x_2,y_1)$ which is quadrirational on
$\CP^1\times\CP^1$. One shows by a direct computation that the maps $F$ for
the above five situations are exactly the five maps listed in Theorem
\ref{Th: intro 1}.

\paragraph{Example of type I.} We use non-homogeneous coordinates $(W_1,W_2)$
on the affine part $\Complex^2$ of $\CP^2$. Consider a pencil of conics 
through four points
$O=(0,0)$, $(0,1)$, $(1,0)$, $(1,1)\in\CP^2$ (any four points on $\CP^2$, no
three of which lie on a straight line, can be brought into these four by a
projective transformation). Conics of this pencil are described by the 
equation (in non-homogeneous coordinates $(W_1,W_2)$ on $\CP^2$):
\[
 Q(\a):\;\; W_2(W_2-1)=\a W_1(W_1-1).
\]
A rational parametrization of such a conic is given, e.g., by
\[
 X(x)=(W_1(x),W_2(x))
     =\Bigl(\frac{x-\a}{x^2-\a},\frac{x(x-\a)}{x^2-\a}\Bigr). 
\]
Here the parameter $x$ has the interpretation of the slope of the line 
$\overline{OX}$. The values of $x$ for the four points of the 
base locus of the pencil on $Q(\a)$ are $x=\a,\infty,0$ and 1. 
A straightforward computation shows that if $Q_1=Q(\a)$, $Q_2=Q(\b)$ then 
the map $F$ defined above coincides with (\ref{F1}).

\paragraph{Example of type V.} Consider a pencil of conics having 
a triple tangency point at $(W_1:W_2:W_3)=(0:1:0)$  
(in homogeneous coordinates on $\CP^2$). Conics of
this pencil and their rational parametrization are given by  formulae
\[
 Q(\a):\;\;W_2-W_1^2-\a=0,\qquad
 X(x)=(W_1(x),W_2(x))=(x,x^2+\a). 
\]
If $Q_1=Q(\a)$, $Q_2=Q(\b)$ then the map $F$ defined above is (\ref{F5}).
\medskip

The following theorem delivers a geometric interpretation of the statement
of Theorem \ref{Th: intro 2}.

\begin{theorem}\label{th:greeks}
Let $Q_i$, $i=1,2,3$ be three non-degenerate members of a linear pencil of
conics. Let $X\in Q_1$, $Y\in Q_2$ and $Z\in Q_3$ be arbitrary
points on these conics. Define the maps $\cF_{ij}$ as above,
corresponding to the pair of conics $(Q_i,Q_j)$. Set
$(X_2,Y_1)=\cF_{12}(X,Y)$, $(X_3,Z_1)=\cF_{13}(X,Z)$, and
$(Y_3,Z_2)=\cF_{23}(Y,Z)$. Then
\begin{gather}
\label{geom 1}
 X_{23}=\overline{X_3Y_3}\cap\overline{X_2Z_2}\in Q_1,\quad
 Y_{13}=\overline{X_3Y_3}\cap\overline{Y_1Z_1}\in Q_2, \\
\label{geom 2}
 Z_{12}=\overline{Y_1Z_1}\cap\overline{X_2Z_2}\in Q_3.
\end{gather}
\end{theorem}

\begin{figure}[htbp]\label{fig:greeks3}
\begin{center}\includegraphics[width=10cm]{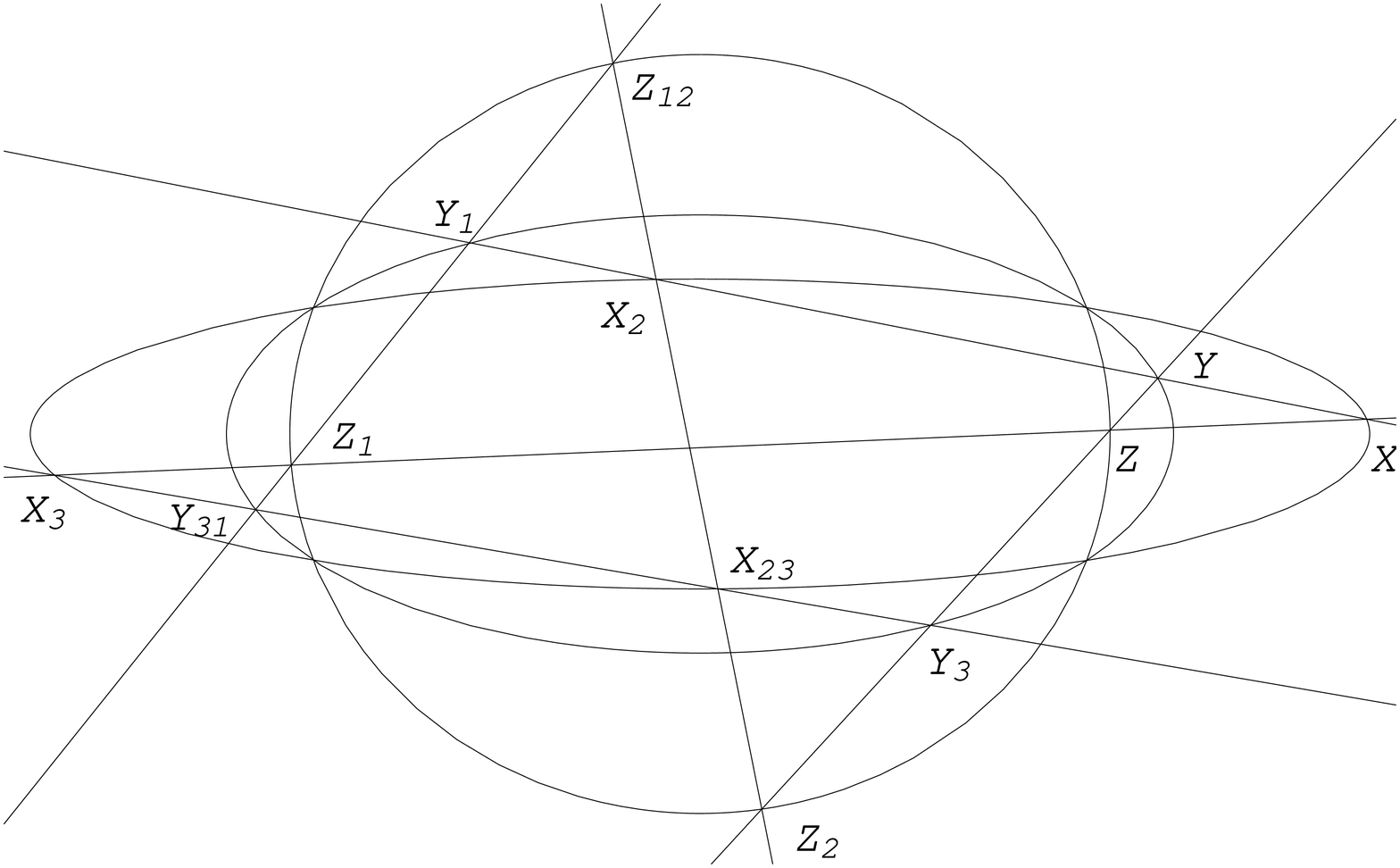}\end{center}
\caption{3D consistency at a linear pencil of conics}
\end{figure}
\begin{proof} 
We will work with equations of lines and conics on $\CP^2$ in homogeneous
coordinates, and use the same notations for geometric objects and homogeneous
polynomials vanishing on these objects. We start with the following data: two
conics $Q_1$ and $Q_2$, three points $X\in Q_1$, $Y\in Q_2$, and $Z$. Construct
the lines $a=\overline{YZ}$, $b=\overline{XZ}$ and  $c=\overline{XY}$,
respectively. Let 
\begin{gather*}
 X_2=(c\cap Q_1)\setminus X,\quad Y_1=(c\cap Q_2)\setminus Y, \\
 X_3=(b\cap Q_1)\setminus X,\quad Y_2=(a\cap Q_2)\setminus Y. 
\end{gather*}
Next, construct the line $C=\overline{X_3Y_3}$, and let 
\[
 X_{23}=(C\cap Q_1)\setminus X_3,\quad Y_{13}=(C\cap Q_2)\setminus Y_3.
\]
Finally, construct the lines $A=\overline{Y_1Y_{13}}$ and
$B=\overline{X_2X_{23}}$.

We have four points $X$, $X_2$, $X_3$ and $X_{23}$ on the conic $Q_1$,
and two pairs of lines $(C,c)$ and $(B,b)$ through two pairs of these points
each. Therefore, there exists $\mu_1\in\CP^1$ such that the conic $Q_1$
has the equation $Q_1=0$ with
\[
 Q_1=\mu_1 Bb+Cc.
\]
Similarly, the conic $Q_2$ has the equation $Q_2=0$ with
\[
 Q_2=\mu_2 Aa+Cc.
\]
Consider the conic
\[
 Q_1-Q_2=\mu_1 Bb-\mu_2 Aa=0.
\]
It belongs to a linear pencil of conics spanned by $Q_1$ and $Q_2$. Further,
the point $Z=a\cap b$ lies on this conic.  Therefore, it has to coincide with
$Q_3$, which has therefore the equation $Q_3=0$ with
\[
 Q_3=\mu_1 Bb-\mu_2 Aa.
\]
Further, the two points $Z_2=a\cap B$ and $Z_1=b\cap A$ also lie on $Q_3$.
Since $Z_2\in B$, we have $B=\overline{X_2Z_2}$. Similarly, since
$Z_1\in A$, we have $A=\overline{Y_1Z_1}$. Finally, we find that the point
$Z_{12}=A\cap B=\overline{Y_1Z_1}\cap\overline{X_2Z_2}$ also lies on $Q_3$, 
which is equivalent to (\ref{geom 2}). 
\end{proof}

\section{Quadrirational and bi-M\"obius maps}\label{s:4rat}

We proceed with a general study of quadrirational maps in the case $\cX=\CP^1$.
In this case we can immediately make some conclusions about their nature. First
of all, for a fixed value of $y$ (with a possible exception of finitely many
values), the maps $F$, $\F$ induce a birational correspondence between
$\CP^1\ni x$ and $\CP^1\ni u$. Similarly, for a fixed value of $x$ (with a
possible exception of finitely many values), the maps $F$, $\F^{-1}$ induce a
birational correspondence between $\CP^1\ni y$ and $\CP^1\ni v$. Therefore, a
quadrirational map has to be {\it bi-M\"obius}, i.e. be given by the formulas
(\ref{map}) from the Introduction. Equivalently, one can say that the variety
$\Gamma_F$ is given by the following two (irreducible) polynomial equations:

\begin{equation}\label{phi eqs}
  F:\;\; \phi(x,y,u)=0\,,\quad \Phi(y,x,v)=0\,,
\end{equation}
where
\begin{align*}
 \phi(x,y,u) &= c(y)xu+d(y)u-a(y)x-b(y), \\
 \Phi(y,x,v) &= C(x)yv+D(x)v-A(x)y-B(x). 
\end{align*}
In order for a bi-M\"obius map to be quadrirational, the following
{\it nondegeneracy conditions} have to be satisfied:
\[
  \text{(i)}\quad u_xv_y-u_yv_x\not\equiv0, 
  \quad\text{(ii)}\quad u_x\not\equiv0,\ v_y\not\equiv0, 
  \quad\text{(iii)}\quad u_y\not\equiv0,\ v_x\not\equiv 0.
\]
Condition (i) is necessary for the existence of the inverse map; conditions 
(ii) are necessary for the existence of the companion maps and are equivalent 
to
\[
 r(y)=a(y)d(y)-b(y)c(y)\not\equiv0,\quad 
 R(x)=A(x)D(x)-B(x)C(x)\not\equiv0;
\]
finally, conditions (iii) assure that the map is really two-dimensional:
if one of them is violated, then, making the corresponding M\"obius
transformations in $u$ or/and $v$, we come to the essentially one-dimensional
map with $u=x$, resp. $v=y$.

The above reasoning holds for the inverse and companion maps as well, so that
for a quadrirational map all four maps $F,F^{-1},\F,\F^{-1}$ have to be
bi-M\"obius. So, the companion maps have to be given by
\begin{alignat}{3}
\label{conj1}
 \F:&\quad & 
    x&= \frac{d(y)u-b(y)}{-c(y)u+a(y)}\,, &\quad
    v&= \frac{\hat{A}(u)y+\hat{B}(u)}{\hat{C}(u)y+\hat{D}(u)}\,, \\
\label{conj2}	
 \F^{-1}:& & 
    y&= \frac{D(x)v-B(x)}{-C(x)v+A(x)}\,, &\quad
    u&= \frac{\hat{a}(v)x+\hat{b}(v)}{\hat{c}(v)x+\hat{d}(v)}\,. 
\end{alignat}
In other words, the variety $\Gamma_F$ can be alternatively described by
either pair of polynomial equations
\begin{alignat*}{3}
 \F&:\quad & \phi(x,y,u) &=0,& \quad \hat\Phi(y,u,v)&=0, \\
 \F^{-1}&: & \Phi(y,x,v) &=0,& \quad \hat\phi(x,v,u)&=0, 
\end{alignat*}
where
\begin{align*}
  \hat\Phi(y,u,v) &= \hat{C}(u)yv+\hat{D}(u)v-\hat{A}(u)y-\hat{B}(u), \\
  \hat\phi(x,v,u) &= \hat{c}(v)xu+\hat{d}(v)u-\hat{a}(v)x-\hat{b}(v).
\end{align*}
We set
\[
 \hat{r}(v)=\hat{a}(v)\hat{d}(v)-\hat{b}(v)\hat{c}(v),\quad
 \hat{R}(u)=\hat{A}(u)\hat{D}(u)-\hat{B}(u)\hat{C}(u).
\]
Finally, the inverse map has to be given by
\[
 F^{-1}:\quad
  x=\frac{\hat{d}(v)u-\hat{b}(v)}{-\hat{c}(v)u+\hat{a}(v)}\,, \quad
  y=\frac{\hat{D}(u)v-\hat{B}(u)}{-\hat{C}(u)v+\hat{A}(u)}\,,
\]
so that the fourth equivalent description of the variety $\Gamma_F$ is
given by the pair of polynomial equations
\[
 F^{-1}:\quad \hat\phi(x,v,u)=0,\quad \hat\Phi(y,u,v)=0.
\]
It is natural to think of the polynomials $\Phi$, $\phi$, $\hat\Phi$ and
$\hat\phi$ as attached to the edges $x,y,u$ and $v$ of the quadrilateral
on Fig.~\ref{fig:maps}, respectively.

For an arbitrary bi-M\"obius map both companion maps are rational. Indeed,
consider, for instance, the companion map $\F:(u,y)\mapsto(x,v)$. The first
equation in (\ref{map}) is immediately solved with respect to $x$, and the
result has to be substituted into the second equation in (\ref{map}), in order
to express $v$ through $u,y$. We come to the following formulas:
\begin{equation}\label{conj1 aux}
 \F:\quad x=X(y,u)=\frac{d(y)u-b(y)}{-c(y)u+a(y)}\,,\quad
         v=\frac{p(y,u)}{q(y,u)}\,.
\end{equation}
This is bi-M\"obius, as in (\ref{conj1}), if and only if the rational
function $p(y,u)/q(y,u)$ is linear-fractional with respect to $y$.
Similarly, the companion map $\F^{-1}:(x,v)\mapsto(u,y)$ is a rational
map of the form
\begin{equation}\label{conj2 aux}
 \F^{-1}:\quad y=Y(x,v)=\frac{D(x)v-B(x)}{-C(x)v+A(x)}\,,\quad
               u=\frac{P(x,v)}{Q(x,v)}\,,
\end{equation}
and is bi-M\"obius, as in (\ref{conj2}), if and only if $P(x,v)/Q(x,v)$ is a
linear-fractional function with respect to $x$. For the inverse map
$(u,v)\mapsto(x,y)$ we have the system
\[
 u=\frac{P(x,v)}{Q(x,v)}\,,\quad v=\frac{p(y,u)}{q(y,u)}\,,
\]
where the first equation has to be solved for $x$, and the second one has to be
solved for $y$. We come to the following conclusion: {\it a bi-M\"obius map
(\ref{map}) is quadrirational, if and only if its both companion maps are
bi-M\"obius. In this case the inverse map is also bi-M\"obius.}

\section{First steps toward classification}\label{s:BD}

According to the nondegeneracy conditions (iii) above, we will
suppose from now on that
\[
 \deg_y\phi\ge1,\quad \deg_x\Phi\ge1.
\]

\begin{lemma}
For a bi-M\"obius map to be quadrirational it is necessary that
\[
 \deg_y\phi\le2,\quad \deg_x\Phi\le2.
\]
\end{lemma}
\begin{proof} 
Eliminating $v$ between $\Phi(y,x,v)=0$ and $\hat\Phi(y,u,v)=0$, we find the
polynomial equation which holds on $\Gamma_F$:
\[
 \Phi\hat\Phi_v-\Phi_v\hat\Phi=0.
\]
The polynomial on the left-hand side depends on $(x,y,u)$ only; it does not
vanish identically (this would violate the nondegeneracy conditions). 
Therefore, the above equation has to be a consequence of  the first equation in
(\ref{map}) (or in (\ref{phi eqs})). In other words, we obtain:
\begin{equation}\label{mu}
 \Phi\hat\Phi_v-\Phi_v\hat\Phi=\mu(x,y,u)\phi,
\end{equation}
with some polynomial factor $\mu\not\equiv0$. The left-hand side
is at most quadratic in $y$. Therefore, $\deg_y\phi\le2$. 
\end{proof}

\begin{lemma}\label{lemma semilinear}
If for a quadrirational map (\ref{map}) $\deg_x\Phi=1$, then
$\deg_u\hat\Phi=1$. Conversely, two arbitrary polynomials $\Phi(y,x,v)$ and
$\hat\Phi(y,u,v)$, each linear in all three arguments, define a quadrirational
map $F:(x,y)\mapsto(u,v)$.
\end{lemma}
\begin{proof} 
Under the condition $\deg_x\Phi=1$, we have for the
polynomials $p(y,u)$, $q(y,u)$ in (\ref{conj1 aux}) the following expressions:
\begin{gather*}
 p(y,u)=(-c(y)u+a(y))(A(X)y+B(X)),\\
 q(y,u)=(-c(y)u+a(y))(C(X)y+D(X)),
\end{gather*}
where $X=(d(y)u-b(y))/(-c(y)u+a(y))$. They are at most linear in $u$, therefore
the same holds for $\hat{A}(u)y+\hat{B}(u)$, $\hat{C}(u)y+\hat{D}(u)$.
Conversely, for any two polynomials $\Phi(y,x,v)$ and $\hat\Phi(y,u,v)$
linear in all arguments, the left--hand side of (\ref{mu}) is linear in $x,u$
and at most quadratic in $y$. If it is irreducible, then it gives
$\phi(x,y,u)$, otherwise it coincides with $\phi$ up to a factor $\mu=\mu(y)$.
This shows that the map $F$ is bi-M\"obius, and the proof for $\F^{\pm1}$ and
$F^{-1}$ is analogous.  
\end{proof}

\paragraph{Remark 1.} It can happen that on the last step in the previous
proof one gets $\phi=\phi(x,u)$, i.e. the dependence on $y$ drops out,
thus leading to degenerate maps $F$ with $u_y=0$ (excluded from
consideration). This happens, for instance, if $\hat\Phi(y,u,v)=\Phi(y,u,v)$,
when $\phi=x-u$.
\medskip

According to Lemma \ref{lemma semilinear}, all quadrirational maps belong
to one of the three subclasses:
\begin{align*}
 \text{[1:1]}&\quad \deg_x\Phi=\deg_u\hat\Phi=\deg_y\phi=\deg_v\hat\phi=1;\\
 \text{[1:2]}&\quad \deg_x\Phi=\deg_u\hat\Phi=1,\quad 
                    \deg_y\phi=\deg_v\hat\phi=2,\quad \text{or vice versa};\\
 \text{[2:2]}&\quad \deg_x\Phi=\deg_u\hat\Phi=\deg_y\phi=\deg_v\hat\phi=2.
\end{align*}
Moreover, Lemma \ref{lemma semilinear} gives an explicit description of
a family of quadrirational maps containing all maps of subclasses [1:1] and 
[1:2], in terms of two polynomials attached to the {\it opposite} (rather than
adjacent) edges of the quadrilateral on Fig.~\ref{fig:maps}. Of course, it 
remains to classify all such maps modulo independent M\"obius transformations 
on $x,y,u,v$. We return to the (simpler) classes [1:1] and [1:2]
in Appendix \ref{a:12}, and concentrate in the main text on the
most interesting subclass [2:2]. So, from now on we assume that
\begin{equation}\label{main assump}
 \deg_x\Phi=2,\quad \deg_y\phi=2.
\end{equation}
It turns out that in this case it is also possible to find an explicit
description of all quadrirational maps in terms of two polynomials
attached to the opposite edges.

\begin{theorem}\label{Th quadratic}
For a quadrirational map (\ref{map}), let $\deg_x\Phi=\deg_y\phi=2$. Then the
polynomials $\Phi(y,x,v)$ and $\hat\Phi(y,u,v)$ corresponding to the
opposite edges of the elementary quadrilateral are related by the formula
\begin{equation}\label{BD2}
 \hat\Phi(y,u,v)=(\g u+\d)^2\Phi\Bigl(y,\frac{\a u+\b}{\g u+\d},v\Bigr).
\end{equation}
Conversely, for an arbitrary polynomial $\Phi(y,x,v)$, linear in $y,v$
and quadratic in $x$, define the polynomial $\hat\Phi(y,u,v)$ as in
(\ref{BD2}). Then the pair of polynomials $(\Phi,\hat\Phi)$ determines
a quadrirational map.
\end{theorem}
\begin{proof} 
Under the assumption $\deg_y\phi=2$, the factor $\mu$ in the formula (\ref{mu})
does not depend on $y$, thus being a polynomial $\mu(x,u)$ of bidegree (1,1).
Moreover, it is irreducible. Indeed, suppose that $\mu(x_0,u)\equiv0$. If
$x_0\neq\infty$, then set $\Psi=\Psi(y,v)=\Phi(y,x_0,v)$, otherwise take
$\Psi=\Psi(y,v)$ equal to the $x^2$-coefficient of $\Phi$. Then Eq. (\ref{mu})
implies $\Psi\hat\Phi_v-\Psi_v\hat\Phi\equiv0$, and by further differentiation
$\Psi\hat\Phi_{uv}-\Psi_v\hat\Phi_u\equiv0$. Therefore, we find:
$\hat\Phi\hat\Phi_{uv}-\hat\Phi_u\hat\Phi_v\equiv0$ --- a contradiction to
the nondegeneracy assumption. So, $\mu$ is of the form
\[
 \mu=\g xu+\d x-\a u-\b,\quad \a\d-\b\g\ne0.
\]
Consider the polynomial
\[
 \bar\Phi(y,u,v)=(\g u+\d)^2\Phi\Big(y,\frac{\a u+\b}{\g u+\d},v\Big).
\]
Eq. (\ref{mu}) implies: $\bar\Phi\hat\Phi_v-\bar\Phi_v\hat\Phi=0$,
therefore $\hat\Phi=\nu(y,u)\bar\Phi$, and since both $\bar\Phi$ and
$\hat\Phi$ are linear in $y$ and quadratic in $u$, we conclude that
$\hat\Phi=\const\bar\Phi$. This proves the first claim of the theorem.

To prove the second one, we start with an arbitrary $\Phi(y,x,v)$ and construct
$\hat\Phi(y,u,v)$ as in (\ref{BD2}). Without loss of generality we can take
simply $\hat\Phi(y,u,v)=\Phi(y,u,v)$ (otherwise applying a M\"obius
transformation in $u$). Then the polynomial $\Phi\hat\Phi_v-\Phi_v\hat\Phi$ is
divisible by $x-u$. In other words, there holds (\ref{mu}) with $\mu=x-u$, and
with the quotient $\phi$ linear in $x,u$ and at most quadratic in $y$. This
demonstrates that the map $F$ is bi-M\"obius, and for  $\F^{\pm1}$ and $F^{-1}$
everything is similar. 
\end{proof}

\paragraph{Remark 2.} 
On the last step of the previous proof, it could well happen that $\phi$ is
less then quadratic in $y$, so that the resulting map is actually in the 
subclass [1:2] ($\deg_x\Phi=\deg_u\hat\Phi=2$ and 
$\deg_y\phi=\deg_v\hat\phi=1$).

\paragraph{Remark 3.} 
A naive counting suggests that the families of quadrirational maps constructed
in this section depends on 2 essential parameters. Indeed, two generic
polynomials $\Phi$, $\hat\Phi$ from the Lemma \ref{lemma semilinear} contain
$2\times8=16$ parameters, 2 are scaled out by homogeneity, and $4\times3=12$
ones correspond to the M\"obius freedom of all variables $x,y,u,v$.
Analogously, the polynomial $\Phi$ from Theorem \ref{Th quadratic} contains 12
parameters, one is scaled out, and $3\times3=9$ ones correspond to the M\"obius
freedom of $x,y,v$.

\section{Singularities}\label{s:sing}

We proceed with the in-depth study of quadrirational maps with the property
(\ref{main assump}), which will eventually lead to their complete
classification. In performing such a classification, it is natural
to factor out the action by the automorphisms group of $\cX$ on all four
variables (edges) $x,y,u,v$ independently. In the present case $\cX=\CP^1$, we
factor out the action by independent M\"obius transformation on all of
$x,y,u,v$. 

According to the conclusion of Sect.~\ref{s:4rat}, we have to study the 
condition under which the companion maps of the given one are bi-M\"obius. 
Recall Eqs. (\ref{conj1 aux}) for the companion map $\F$. Under the
condition (\ref{main assump}), the polynomials $p(y,u)$, $q(y,u)$
are given by
\begin{equation}\label{conj aux}
\begin{aligned}
 p(y,u)&=(-c(y)u+a(y))^2(A(X)y+B(X)), \\
 q(y,u)&=(-c(y)u+a(y))^2(C(X)y+D(X)). 
\end{aligned}
\end{equation}
They are of degree $5$ in $y$ and of degree $2$ in $u$. Similarly,
in Eqs. (\ref{conj2 aux}) for the companion map $\F^{-1}$ the polynomials
$P(x,v)$, $Q(x,v)$ are of degree $5$ in $x$ and of degree $2$ in $v$.

Obviously, in order that the companion map $\F:(y,u)\mapsto(v,x)$ be
bi-M\"obius of the form (\ref{conj1}), it is necessary and sufficient that
\begin{equation}\label{conj1 factor}
 p(y,u)=\rho(y)(\hat{A}(u)y+\hat{B}(u)),\quad
 q(y,u)=\rho(y)(\hat{C}(u)y+\hat{D}(u)),
\end{equation}
where $\rho(y)$ is a polynomial of degree 4. The polynomials $\hat{A}(u),
\ldots,\hat{D}(u)$ are quadratic. Similarly, for the companion map
$\F^{-1}:(v,x)\mapsto(y,u)$ to be bi-M\"obius of the form (\ref{conj2}), it is
necessary and sufficient that
\begin{equation}\label{conj2 factor}
 P(x,v)=\Rho(x)(\hat{a}(v)x+\hat{b}(v)),\quad
 Q(x,v)=\Rho(x)(\hat{c}(v)x+\hat{d}(v)),
\end{equation}
where $\Rho(x)$ is a polynomial of degree 4, and $\hat{a}(v),\ldots,\hat{d}(v)$
are quadratic polynomials.

Next, consider the conditions (\ref{conj1 factor}), (\ref{conj2 factor}). It is
enough to consider the first one, since for the second everything is similar.
For (\ref{conj1 factor}) to hold, there have to exist $4$ numbers $y_i$,
$i=1,\ldots,4,$ such that $p(y_i,u)\equiv0$, $q(y_i,u)\equiv0$ (of course, the
zeros $y_i$ have to be counted with their multiplicity).

\begin{lemma}\label{lemma conj bimob}
If $p(\eta,u)\equiv0$, $q(\eta,u)\equiv0$ for some $\eta\in\CP^1$, then $\eta$
is with necessity a root of $r(y)=a(y)d(y)-b(y)c(y)$.
\end{lemma}
\begin{proof} 
From the formulas (\ref{conj aux}) one sees that if
$a(\eta)d(\eta)-b(\eta)c(\eta)\ne0$ then
\begin{align*}
 A\left(\frac{du-b}{-cu+a}\right)\eta+
 B\left(\frac{du-b}{-cu+a}\right) & \equiv0, \\
 C\left(\frac{du-b}{-cu+a}\right)\eta+
 D\left(\frac{du-b}{-cu+a}\right) & \equiv0,
\end{align*}
where $a=a(\eta),\ldots,d=d(\eta)$. These two formulas, being identities in
$u$, are equivalent to $A(x)\eta+B(x)\equiv0$, $C(x)\eta+D(x)\equiv0$, which,
in turn, means that the second component of the map (\ref{map}) is degenerate
($R(x)\equiv0$) --- a contradiction. Hence, $a(\eta)d(\eta)-b(\eta)c(\eta)=0$.
Lemma is proved.  
\end{proof}

\begin{proposition}[cancellation of simple zero]\label{simple root}
In order for $\eta$ to be a common zero of
$p(y,u)$ and $q(y,u)$ (considered as polynomials in $y$), i.e.
\[
 p(\eta,u)=q(\eta,u)=0,
\]
it is necessary and sufficient that the point $(\xi,\eta)$ defined by the
system
\begin{equation}\label{system1}
 a(\eta)\xi+b(\eta)=0,\quad c(\eta)\xi+d(\eta)=0,
\end{equation}
satisfies also the system
\begin{equation}\label{system2}
 A(\xi)\eta+B(\xi)=0,\quad C(\xi)\eta+D(\xi)=0.
\end{equation}
In particular, $\xi$ is a zero of $R(x)$, and moreover a common zero of
$P(x,v)$, $Q(x,v)$, considered as polynomials of $x$:
\[
 P(\xi,v)=Q(\xi,v)=0.
\]
\end{proposition}
\begin{proof} 
Since $\eta$ has to be a root of $r(y)$, i.e. $r(\eta)=0$, we see that the
system (\ref{system1}) admits a nontrivial solution. (Of course, in this system
the point $(\xi,1)$ is considered as an element of $\CP^1$, so that in the case
$a(\eta)=c(\eta)=0$ it has to be taken equal to $(1,0)$.) First suppose that
$(a(\eta),c(\eta))\ne(0,0)$. Notice that there holds, identically in $u$,
\begin{equation}\label{degen}
 X(\eta,u)=\frac{d(\eta)u-b(\eta)}{-c(\eta)u+a(\eta)}\equiv \xi\,.
\end{equation}
Substituting this into (\ref{conj aux}) we find:
\begin{align*}
 p(\eta,u) & = (-c(\eta)u+a(\eta))^2(A(\xi)\eta+B(\xi))\,, \\
 q(\eta,u) & = (-c(\eta)u+a(\eta))^2(C(\xi)\eta+D(\xi))\,, 
\end{align*}
and these both polynomials in $u$ vanish identically, if and only
if (\ref{system2}) holds. Now let $a(\eta)=c(\eta)=0$. This
means that the corresponding $\xi=\infty$, so that the system
(\ref{system2}) has to be read projectively, as:
\begin{equation}\label{system2 mod}
 \tilde{A}(0)\eta+\tilde{B}(0)=0\,,\quad
 \tilde{C}(0)\eta+\tilde{D}(0)=0\,,
\end{equation}
where $\tilde{A}(x)=x^2A(1/x)$, etc. At the same time, the formulas 
(\ref{conj aux}) have to be replaced by
\begin{align*}
 p(y,u) &= (d(y)u-b(y))^2\bigl(\tilde{A}(1/X)y+\tilde{B}(1/X)\bigr),\\
 q(y,u) &= (d(y)u-b(y))^2\bigl(\tilde{C}(1/X)y+\tilde{D}(1/X)\bigr),
\end{align*}
which yield
\begin{align*}
 p(\eta,u) &= (d(\eta)u-b(\eta))^2\bigl(\tilde{A}(0)\eta+\tilde{B}(0)\bigr),\\
 q(\eta,u) &= (d(\eta)u-b(\eta))^2\bigl(\tilde{C}(0)\eta+\tilde{D}(0)\bigr),
\end{align*}
and again these both polynomials in $u$ vanish identically, if and only if
(\ref{system2 mod}) holds. Finally, the statements about $\xi$ follow due to
the symmetry of systems (\ref{system1}), (\ref{system2}) with respect to
changing the roles of $x$ and $y$.  
\end{proof}

Notice that the system (\ref{system1}) characterizes the points $(\xi,\eta)$
where the first fraction in (\ref{map}) is not defined. We see that for
quadri-rational maps these points satisfy with necessity also the system
(\ref{system2}) which characterizes the points where the second fraction in
(\ref{map}) is not defined. Due to Proposition \ref{simple root}, we see that
for quadri-rational maps singularities of both fractions in (\ref{map})
coincide. Clearly, there are no more than four such singularities.

In what follows, we shall take the freedom of performing M\"obius
transformations in the variables $x,y,u,v$ in order not to consider various
particular cases as in the previous proof. Indeed, performing a M\"obius
transformation in $x$, we can achieve that neither root $\xi$ of $R(x)$ lies at
$\infty$, that is, $(a(\eta),c(\eta))\ne(0,0)$ for all zeros $\eta$ of $r(y)$.
Next, by doing a M\"obius transformation in $u$, we can even assure that  both
$a(\eta)\ne0$, $c(\eta)\ne0$. Similarly, we shall assume that both
$A(\xi)\ne0$, $C(\xi)\ne0$. Note also that under these conditions multiple
zeros of the functions $p(y,u)$, $q(y,u)$, considered as polynomials in $y$,
can be alternatively characterized as zeros of the same multiplicity of the
functions
\begin{equation}\label{tilde pq}
 \tilde{p}(y,u)=A(X)y+B(X),\quad \tilde{q}(y,u)=C(X)y+D(X).
\end{equation}
This is seen immediately from (\ref{conj aux}).

\begin{proposition}[cancellation of double zero]\label{double root}
In order for $\eta$ to be a common double zero
of $p(y,u)$ and $q(y,u)$ (considered as polynomials in $y$), 
it is necessary and sufficient that the point $(\xi,\eta)$ satisfies, in
addition to (\ref{system1}), (\ref{system2}), the following equations:
\begin{align}
\label{system3 1}
 &\left\{\begin{array}{l} a'(\eta)\xi+b'(\eta)+\l a(\eta)=0\,,\\
 c'(\eta)\xi+d'(\eta)+\l c(\eta)=0\,, \end{array}\right.\\ 
\label{system3 2}
 &\left\{\begin{array}{l} A'(\xi)\eta+B'(\xi)+\l^{-1}A(\xi)=0\,, \\
 C'(\xi)\eta+D'(\xi)+\l^{-1}C(\xi)=0\,,\end{array}\right.
\end{align}
with some $\l\in\CP^1$. Such $\eta$ is with necessity a double zero of $r(y)$,
and the corresponding $\xi$ is a double zero of $R(x)$. Moreover, $\xi$ is a
common double zero of $P(x,v)$ and $Q(x,v)$ (considered as polynomials in $x$).
\end{proposition}
\begin{proof} 
Differentiate the definitions (\ref{tilde pq}) with respect to $y$, then set
$y=\eta$, taking  into account Eq. (\ref{degen}). The result reads:
\begin{equation}\label{dou aux}
\begin{aligned} 
 (A'(\xi)\eta+B'(\xi))X_y(\eta,u)+A(\xi)&= 0,\\
 (C'(\xi)\eta+D'(\xi))X_y(\eta,u)+C(\xi)&= 0.
\end{aligned}
\end{equation}

Eqs. (\ref{dou aux}) can be fulfilled identically in $u$ if and only if
$X_y(\eta,u)$ is constant, i.e. does not depend on $u$. Denote this constant by
\begin{equation}\label{dou aux4}
 X_y(\eta,u)\equiv\l,
\end{equation}
then Eqs. (\ref{dou aux}) coincide with (\ref{system3 2}). Differentiate
\begin{equation}\label{X aux}
 X(y,u)(c(y)u-a(y))+d(y)u-b(y)=0
\end{equation}
with respect to $y$ and set $y=\eta$. Taking into account that
$X(\eta,u)\equiv\xi$ and $X_y(\eta,u)\equiv\l$, one comes to the formula
\[
 \l(cu-a)+\xi(c'u-a')+(d'u-b')=0,
\]
where we write for shortness $a$, $a'$ etc. for $a(\eta)$, $a'(\eta)$ etc.
This is equivalent to (\ref{system3 1}). To prove
that $\eta$ is a double root of $r(y)$, calculate:
\[
 r'(\eta)=a'd+ad'-b'c-bc'=-(a'\xi+b')c+a(c'\xi+d')=0;
\]
on the last step we used Eq. (\ref{system3 1}). The proof that $\xi$ is
a double root of $R(x)$ is analogous, since the systems (\ref{system3 1}),
(\ref{system3 2}) are symmetric with respect to the interchange
$\xi\leftrightarrow\eta$. 
\end{proof}

It will be convenient to represent the parameter $\l$ as a quotient:
\begin{equation}\label{lambda}
 \l=\dot\xi/\dot\eta.
\end{equation}
Clearly, the quantities $\dot\xi$, $\dot\eta$ are defined only
up to a common factor. They can be considered as representing the velocity
vector at $(\xi,\eta)$ of a parametrized curve on the plane
$\CP^1\times\CP^1$ passing through this point. The freedom of a
simultaneous multiplying $\dot\xi$ and $\dot\eta$ by a common factor
corresponds to a possibility to re-parametrize such a curve.
This supports the intuitive understanding of the present case
as a degeneration of the previous one, when two ``simple singularities''
glue together moving along some curve with a certain limiting velocity.
This understanding is further supported by the form of the system
(\ref{system3 1}), (\ref{system3 2}) in this new notation:
\begin{align}\label{system31 new}
 &\left\{\begin{array}{l}
 (a'(\eta)\xi+b'(\eta))\dot\eta+a(\eta)\dot\xi=0,\\
 (c'(\eta)\xi+d'(\eta))\dot\eta+c(\eta)\dot\xi=0,
 \end{array}\right.\\ 
 &\left\{\begin{array}{l}
 (A'(\xi)\eta+B'(\xi))\dot\xi+A(\xi)\dot\eta=0,\\
 (C'(\xi)\eta+D'(\xi))\dot\xi+C(\xi)\dot\eta=0.
 \end{array}\right.
\end{align}
This can be clearly read as the result of differentiation of
(\ref{system1}), (\ref{system2}) at $(\xi,\eta)$ along the parametrized
curve through $(\xi,\eta)$ with the velocity vector $(\dot\xi,\dot\eta)$.

\begin{proposition}[cancellation of triple zero]\label{triple root}
In order for $\eta$ to be a common triple zero of $p(y,u)$ and $q(y,u)$
(considered as polynomials in $y$), it is necessary and sufficient that the
point $(\xi,\eta)$ satisfies, in addition to (\ref{system1}), (\ref{system2}),
(\ref{system3 1}), (\ref{system3 2}), the following equations:
\begin{align}
\label{system4 1}
 &\left\{\begin{array}{ll}
   a''(\eta)\xi+b''(\eta)+2\l a'(\eta)+\l^{3/2}\theta a(\eta)&=0,\\
   c''(\eta)\xi+d''(\eta)+2\l c'(\eta)+\l^{3/2}\theta c(\eta)&=0,
  \end{array}\right.\\ 
\label{system4 2}
 &\left\{\begin{array}{ll}
   A''(\xi)\eta+B''(\xi)+2\l^{-1}A'(\xi)-\l^{-3/2}\theta A(\xi)=0,\\
   C''(\xi)\eta+D''(\xi)+2\l^{-1}C'(\xi)-\l^{-3/2}\theta C(\xi)=0,
  \end{array}\right. 
\end{align}
with some $\theta\in\CP^1$. Such $\eta$ is with necessity a triple
zero of $r(y)$, and the corresponding $\xi$ is a triple zero of $R(x)$.
Moreover, $\xi$ is a common triple zero of $P(x,v)$ and $Q(x,v)$
(considered as polynomials in $x$).
\end{proposition}
\begin{proof} 
Differentiate (\ref{tilde pq}) twice with respect to $y$,
then set $y=\eta$, taking into account Eq. (\ref{degen}). The result reads:
\begin{align*}
 & (A''(\xi)\eta+B''(\xi))X_y^2(\eta,u)+2A'(\xi)X_y(\eta,u)\\
 &\qquad\qquad\qquad +(A'(\xi)\eta+B'(\xi))X_{yy}(\eta,u) = 0, \\
 & (C''(\xi)\eta+D''(\xi))X_y^2(\eta,u)+2C'(\xi)X_y(\eta,u)\\
 &\qquad\qquad\qquad +(C'(\xi)\eta+D'(\xi))X_{yy}(\eta,u) = 0.
\end{align*}
Using (\ref{system3 2}) and the value $X_y(\eta,u)=\l$ from (\ref{dou aux4}),
we bring this into the form
\begin{equation}\label{tri aux}
\begin{aligned}
 A''(\xi)\eta+B''(\xi)+2\l^{-1}A'(\xi)-\l^{-3}A(\xi)X_{yy}(\eta,u) &= 0,\\
 C''(\xi)\eta+D''(\xi)+2\l^{-1}C'(\xi)-\l^{-3}C(\xi)X_{yy}(\eta,u) &= 0.
\end{aligned}
\end{equation}
Eqs. (\ref{tri aux}) can be fulfilled identically in $u$ if and only if
$X_{yy}(\eta,u)$ does not depend on $u$; we denote this constant by
\begin{equation}\label{tri aux4}
X_{yy}(\eta,u)\equiv \l^{3/2}\theta\,,
\end{equation}
and then (\ref{tri aux}) coincide with (\ref{system4 2}). Next, differentiate
(\ref{X aux}) twice with respect to $y$ and set $y=\eta$. Taking into account
Eqs. (\ref{degen}), (\ref{dou aux4}) and (\ref{tri aux4}), one finds:
\[
 \l^{3/2}\theta(cu-a)+2\l(c'u-a')+\xi(c''u-a'')+(d''u-b'')=0.
\]
This is equivalent to (\ref{system4 1}). It remains to prove that $\eta$
is indeed a triple root of $r(y)$. For this, note that under the
conditions (\ref{system1}), (\ref{system3 1}), we have:
\begin{align*}
 r''(\eta) &= a''d+2a'd'+ad''-b''c-2b'c'-bc''\\
  &= -(a''\xi+b'')c+(c''\xi+d'')a+2a'(-c'\xi-\l c)-2c'(-a'\xi-\l a)\\
  &= -(a''\xi+b''+2\l a')c+(c''\xi+d''+2\l c')a,
\end{align*}
and this vanishes due to (\ref{system4 1}).  
\end{proof}

It will be convenient to parametrize $\theta$ by the 2--germ of a curve through
$(\xi,\eta)$, similar to the parametrization (\ref{lambda}) of $\l$ by the
1--germ. To find such a parametrization, differentiate the relations
(\ref{system1}), (\ref{system2}) twice along a (parametrized) curve through
$(\xi,\eta)$ with the 2--germ at this point given by
$(\dot\xi,\ddot\xi,\dot\eta,\ddot\eta)$:
\begin{align}\label{system41 new}
 &\left\{\begin{array}{l}
   (a''(\eta)\xi+b''(\eta))\dot\eta^2+2a'(\eta)\dot\xi\dot\eta+
   (a'(\eta)\xi+b'(\eta))\ddot\eta+a(\eta)\ddot\xi=0,\\
   (c''(\eta)\xi+d''(\eta))\dot\eta^2+2c'(\eta)\dot\xi\dot\eta+
   (c'(\eta)\xi+d'(\eta))\ddot\eta+c(\eta)\ddot\xi=0,
  \end{array}\right. \\
 &\left\{\begin{array}{l}
  (A''(\xi)\eta+B''(\xi))\dot\xi^2+2A'(\xi)\dot\xi\dot\eta+
  (A'(\xi)\eta+B'(\xi))\ddot\xi+A(\xi)\ddot\eta=0,\\
  (C''(\xi)\eta+D''(\xi))\dot\xi^2+2C'(\xi)\dot\xi\dot\eta+
  (C'(\xi)\eta+D'(\xi))\ddot\xi+C(\xi)\ddot\eta=0.
  \end{array}\right. 
\end{align}
Taking into account Eqs. (\ref{system31 new}) and (\ref{lambda}), one sees that
the above relations can be written as (\ref{system4 1}), (\ref{system4 2}),
respectively, provided
\begin{equation}\label{theta}
 \theta=\frac{\ddot\xi\dot\eta-\dot\xi\ddot\eta}
             {(\dot\xi\dot\eta)^{3/2}}\,.
\end{equation}
This representation is actually independent of the parametrization
of the curve, and respects the symmetry $\theta\to -\theta$ as
$\xi\leftrightarrow\eta$, as it should.

\begin{proposition}[cancellation of quadruple zero]\label{quadruple root}
In order for $\eta$ to be a common quadruple zero of $p(y,u)$ and $q(y,u)$
(considered as polynomials in $y$),  it is necessary and sufficient that the
point $(\xi,\eta)$ satisfies, in addition to (\ref{system1}), (\ref{system2}),
(\ref{system3 1}), (\ref{system3 2}), (\ref{system4 1}), (\ref{system4 2}), the
following equations:
\begin{align}
\label{system5 1}
 &\left\{\begin{array}{l}
  a''(\eta)+\l^{1/2}\theta a'(\eta)+\frac12\l(\theta^2+\k)a(\eta)=0,\\
  c''(\eta)+\l^{1/2}\theta c'(\eta)+\frac12\l(\theta^2+\k)c(\eta)=0,
  \end{array}\right.\\
\label{system5 2}
 &\left\{\begin{array}{l}
  A''(\xi)-\l^{-1/2}\theta A'(\xi)+\frac12\l^{-1}(\theta^2-\k)A(\xi)=0,\\
  C''(\xi)-\l^{-1/2}\theta C'(\xi)+\frac12\l^{-1}(\theta^2-\k)C(\xi)=0,
\end{array}\right. 
\end{align}
with some $\k\in\CP^1$. Such $\eta$ is with necessity a quadruple
zero of $r(y)$, and the corresponding $\xi$ is a quadruple zero of $R(x)$.
Moreover, $\xi$ is a common quadruple zero of $P(x,v)$ and $Q(x,v)$
(considered as polynomials in $x$).
\end{proposition}
\begin{proof} 
Differentiate (\ref{tilde pq}) thrice with respect to $y$, then set $y=\eta$,
taking into account Eq. (\ref{degen}). The result reads:
\begin{align*}
 (A''(\xi)\eta+B''(\xi))X_y(\eta,u)X_{yy}(\eta,u)+A'(\xi)X_{yy}(\eta,u)
  +A''(\xi)X_y^2(\eta,u) & \\
 +\tfrac13(A'(\xi)\eta+B'(\xi))X_{yyy}(\eta,u) &= 0, \\
 (C''(\xi)\eta+D''(\xi))X_y(\eta,u)X_{yy}(\eta,u)+C'(\xi)X_{yy}(\eta,u)
  +C''(\xi)X_y^2(\eta,u) & \\
 +\tfrac13(C'(\xi)\eta+D'(\xi))X_{yyy}(\eta,u) &= 0.
\end{align*}
Using (\ref{system3 2}), (\ref{system4 2}), and the values
$X_y(\eta,u)=\l$, $X_{yy}(\eta,u)=\l^{3/2}\theta$ from
(\ref{dou aux4}), (\ref{tri aux4}), we bring this into the form
\begin{equation}\label{qua aux}
\begin{aligned}
  A''(\xi)-\l^{-1/2}\theta A'(\xi)+\l^{-1}\theta^2 A(\xi)
     -\tfrac13\l^{-3}A(\xi)X_{yyy}(\eta,u) &= 0, \\
  C''(\xi)-\l^{-1/2}\theta C'(\xi)+\l^{-1}\theta^2 C(\xi)
     -\tfrac13\l^{-3}C(\xi)X_{yyy}(\eta,u) &= 0.
\end{aligned}
\end{equation}
Eqs. (\ref{qua aux}) are fulfilled identically in $u$ if and only if
$X_{yyy}(\eta,u)$ does not depend on $u$. Denote this constant by
\begin{equation}\label{qua aux4}
 X_{yyy}(\eta,u)\equiv \frac{3\l^2}{2}(\theta^2+\k)\,,
\end{equation}
Then (\ref{qua aux}) are equivalent to (\ref{system5 2}). Next, differentiate
(\ref{X aux}) thrice with respect to $y$ and set $y=\eta$. Taking into account
(\ref{degen}), (\ref{dou aux4}), (\ref{tri aux4}) and (\ref{qua aux4}), we
find:
\[
 \tfrac12\l^2(\theta^2+\k)(cu-a)+\l^{3/2}\theta(c'u-a')+\l(c''u-a'')=0,
\]
which is equivalent to (\ref{system5 1}).
It remains to prove that $\eta$ is indeed a quadruple root of $r(y)$.
For this, note:
\begin{align*}
 \tfrac13r'''(\eta) 
  &= a''d'+a'd''-b''c'-b'c''\\
  &= a''d'-a'(c''\xi+2\l c'+\l^{3/2}\theta c)
     +(a''\xi+2\l a'+\l^{3/2}\theta a)c'-b'c''\\
  &= a''(c'\xi+d')-c''(a'\xi+b')+\l^{3/2}\theta(ac'-a'c)\\
  &= -\l c(a''+\l^{1/2}\theta a')+\l a(c''+\l^{1/2}\theta c')=0.
\end{align*}
In this computation we used subsequently (\ref{system4 1}), (\ref{system3 1}),
(\ref{system5 1}). 
\end{proof}

Like in the cases of double and triple zeros (cf. (\ref{lambda}),
(\ref{theta})),
it will be convenient to parametrize $\k$ by the 3--germ of
a curve through $(\xi,\eta)$. To find the corresponding formula,
differentiate the relations (\ref{system1}), (\ref{system2})
thrice along a (parametrized) curve through $(\xi,\eta)$
with the 3--germ at this point given by
$(\dot\xi,\ddot\xi,\dddot{\xi},\dot\eta,\ddot\eta,\dddot{\eta})$:
\begin{align*}
 &\left\{\begin{array}{l}
   (a''(\eta)\xi+b''(\eta))\ddot\eta\dot\eta
   +a''(\eta)\dot\xi\dot\eta^2+a'(\eta)(\ddot\xi\dot\eta+\dot\xi\ddot\eta)\\
  \qquad\qquad\qquad
   +\frac13\Bigl((a'(\eta)\xi+b'(\eta))\dddot\eta+a(\eta)\dddot\xi\Bigr)=0,\\
   (c''(\eta)\xi+d''(\eta))\ddot\eta\dot\eta
   +c''(\eta)\dot\xi\dot\eta^2+c'(\eta)(\ddot\xi\dot\eta+\dot\xi\ddot\eta)\\
  \qquad\qquad\qquad
   +\frac13\Bigl((c'(\eta)\xi+d'(\eta))\dddot\eta+c(\eta)\dddot\xi\Bigr)=0\,
  \end{array}\right. \\
 &\left\{\begin{array}{l}
   (A''(\xi)\eta+B''(\xi))\ddot\xi\dot\xi
   +A''(\xi)\dot\xi^2\dot\eta+A'(\xi)(\ddot\xi\dot\eta+\dot\xi\ddot\eta)\\
  \qquad\qquad\qquad
   +\frac13\Bigl((A'(\xi)\eta+B'(\xi))\dddot\xi+A(\xi)\dddot\eta\Big)=0,\\
   (C''(\xi)\eta+D''(\xi))\ddot\xi\dot\xi
   +C''(\xi)\dot\xi^2\dot\eta+C'(\xi)(\ddot\xi\dot\eta+\dot\xi\ddot\eta)\\
  \qquad\qquad\qquad
   +\frac13\Bigl((C'(\xi)\eta+D'(\xi))\dddot\xi+D(\xi)\dddot\eta\Bigr)=0.
  \end{array}\right. 
\end{align*}
Taking into account Eqs. (\ref{system31 new}), (\ref{system41 new}) and
(\ref{lambda}), (\ref{theta}), one sees that the above relations
can be written as (\ref{system5 1}), (\ref{system5 2}), respectively, provided
\begin{equation}\label{kappa}
 \k= \frac{2(\dddot\xi\dot\eta-\dot\xi\dddot\eta)}
         {3(\dot\xi\dot\eta)^2}
    -\frac{\ddot\xi^2\dot\eta^2-\dot\xi^2\ddot\eta^2}{(\dot\xi\dot\eta)^3}
   =\frac{2}{3\dot\xi\dot\eta}(S(\xi)-S(\eta)),
\end{equation}
where $S(\xi)$ stands for the Schwarzian derivative,
$S(\xi)=\dddot\xi/\dot\xi-3\ddot\xi^2/2\dot\xi^2$. Like (\ref{theta}),
representation (\ref{kappa}) is actually independent of the parametrization of
the curve, and respects the symmetry $\k\to -\k$ as $\xi\leftrightarrow\eta$,
as it should.

\section{Five types of quadrirational maps}\label{s:types}

A polynomial $R(x)$ of the 4th degree belongs to one of the five types,
$\cT=\I,\dots,\V$, according to the multiplicity of zeros, as shown in the
second column of the table \ref{tab:inv}. We will denote by $m$ the number of
the distinct roots, given in the third column. The ordering of the roots with
respect to their multiplicity will be firmly associated to the type of the
polynomial. Note that we consider the polynomials projectively, with the action
of the M\"obius transformations given by $R(x)\mapsto (\g x+\d)^4R((\a
x+\b)/(\g x+\d))$. Hence, for instance, the polynomial $x(x-1)$ considered as
the polynomial of the 4th degree is of the type II, having two simple roots
$x_1=0$, $x_2=1$ and one double root $x_3=\infty$.

\def\vv{\vrule height1.8em width0em depth1.2em}
\begin{table}[!ht]
\[
\begin{array}{|c|c|c|c|c|c|}
\hline
 \cT  & \text{zeroes} & m & \cD_\cT(x) & q_\cT(x) & \cC_\cT(\a) \\
\hline
 \I   & (x_1,x_2,x_3,x_4) 
      & 4 
      & (x_1,x_2,x_3,x_4)      
      & \dfrac{(x_1-x_2)(x_3-x_4)}{(x_2-x_3)(x_4-x_1)}
      & (\infty,1,0,\a) \vv\\
 \II  & (x_1,x_2,x_3,x_3) 
      & 3 
      & (x_1,x_2,x_3,\dot x_3) 
      & \dfrac{\dot x_3(x_1-x_2)}{(x_2-x_3)(x_1-x_3)}     
      & (\infty,1,0,\a) \vv\\
 \III & (x_1,x_1,x_2,x_2) 
      & 2 
      & (x_1,\dot x_1,x_2,\dot x_2) 
      & \dfrac{\dot x_1\dot x_2}{(x_1-x_2)^2}
      & (1,1,0,\a) \vv\\
 \IV  & (x_1,x_2,x_2,x_2) 
      & 2 
      & (x_1,x_2,\dot x_2,\ddot x_2) 
      & \dfrac{\ddot x_2}{2\dot x_2}+\dfrac{\dot x_2}{x_1-x_2}
      & (\infty,0,1,2\a) \vv\\  
 \V   & (x_1,x_1,x_1,x_1)  
      & 1 
      & (x_1,\dot x_1,\ddot x_1,\dddot x_1) 
      & \dfrac16\Bigl(\dfrac{\dddot x_1}{\dot x_1}
         -\dfrac{3\ddot x_1^2}{2\dot x_1^2}\Bigr)
      & (0,1,0,6\a) \vv\\
 \hline  
\end{array}
\]
\caption{M\"obius invariants and canonical data}\label{tab:inv}
\end{table}

\begin{theorem}\label{Th sing}
Consider a bi-M\"obius map $F$ defined by (\ref{phi eqs}) with
$\deg_y\phi=\deg_x\Phi=2$. Necessary and sufficient conditions for $F$ to be
quadrirational read:
\begin{itemize}\setlength{\itemsep}{0pt}
\item $R(x)$ and $r(y)$ belong simultaneously to the same type $\cT$;
\item The roots $x_i$, $y_i$ ($i=1,\ldots,m$) of $R(x)$, $r(y)$
      can be ordered so that:
\begin{enumerate}\setlength{\itemsep}{0pt}
\item there hold (\ref{system1}), (\ref{system2}) with 
      $(\xi,\eta)=(x_i,y_i)$ for all $i=1,\ldots,m$;
\item for roots $x_i$, $y_i$ of multiplicity $\ge 2$ there hold
      additionally (\ref{system3 1}), (\ref{system3 2}) with
      $(\xi,\eta)=(x_i,y_i)$ and with some $\l=\l_i\in{\CP}^1$;
\item for roots $x_i$, $y_i$ of multiplicity $\ge 3$ there hold
      additionally (\ref{system4 1}), (\ref{system4 2}) with
      $(\xi,\eta)=(x_i,y_i)$, $\l=\l_i$, and with some
      $\theta=\theta_i\in{\CP}^1$;
\item for quadruple roots $x_i$, $y_i$ there hold additionally 
      (\ref{system5 1}), (\ref{system5 2}) with $(\xi,\eta)=(x_i,y_i)$, 
      $\l=\l_i$, $\theta=\theta_i$ and with some 
      $\k=\k_i\in{\CP}^1$.
\end{enumerate} 
\end{itemize} 
\end{theorem}
\begin{proof}
This follows directly from Propositions \ref{simple root}, \ref{double root},
\ref{triple root}, \ref{quadruple root}: we see that the polynomial $\rho(y)$
in (\ref{conj1 factor}) has to coincide with $r(y)$, while the polynomial
$\Rho(x)$ in (\ref{conj2 factor}) has to coincide with $R(x)$, and, moreover,
both polynomials have to belong {\it simultaneously} to one of the types \I--\V.
\end{proof}

\begin{definition}
We say that a quadrirational map F defined by (\ref{phi eqs}) with
$\deg_y\phi=\deg_x\Phi=2$ is of type $\cT$ if the polynomials $R(x)$ and $r(y)$
are of type $\cT$.
\end{definition}

The formulas
\begin{equation}\label{params}
 \l_i=\frac{\dot x_i}{\dot y_i}\,,\quad
 \theta_i=\frac{\ddot x_i\dot y_i-\dot x_i\ddot y_i}
 {(\ddot x_i\dot y_i)^{3/2}}\,,\quad
 \k_i=\frac{2}{3\dot x_i\dot y_i}(S(x_i)-S(y_i))
\end{equation}
allow us to describe the singularities of a quadrirational map by associating
the data $\cD_\cT(x)$, $\cD_\cT(y)$ to the edges $x,y$ of the elementery
quadrilateral. These data are shown in the fourth column of the Table
\ref{tab:inv},

In all cases except for $\cT=\I$ these data are defined not uniquely, but up to
transformations
\begin{equation}\label{trafo}
 x_i'=\dot x_i\g_1\,,\quad
 x_i''=\ddot x_i\g_1^2+\dot x_i\g_2\,,\quad
 x_i'''=\dddot x_i\g_1^3+3\ddot x_i\g_2\g_1+\dot x_i\g_3
\end{equation}
and similarly for the edge $y$, with the same $\g_1,\g_2,\g_3$. In the only
case $\cT=\III$ when such transformations are applicable with two different
values $i=1,2$, the quantities $\g_1$ can be chosen for these $i$ independently
(but for each $i=1,2$ they still have to be the same for the edges $x$ and
$y$). These quantities $\g_1,\g_2,\g_3$ have the meaning of
\[
 \g_1=t',\quad \g_2=t'',\quad \g_3=t'''
\]
for the change of ``time'' $t=t(s)$ in the parametrized curves
$(x(t),y(t))$ through $(x_i,y_i)$ etc. Clearly, the values of the
parameters (\ref{params}) do not change under transformations
(\ref{trafo}).

The action of M\"obius transformations $x\mapsto(\a x+\b)/(\g x+\d)$ on the
data $\cD_\cT$ is defined by the natural formulas, obtained simply by
differentiation:
\[
 x_i\mapsto\frac{\a x_i+\b}{\g x_i+\d}\,,\quad
 \dot x_i\mapsto\frac{(\a\d-\b\g)\dot x_i}{(\g x_i+\d)^2}\,,\quad
 \text{etc.}
\]

\begin{proposition}
The data $\cD_\cT(x)$ admit the M\"obius invariants $q_\cT(x)$, as shown in the
fifth column of the Table \ref{tab:inv}.
\end{proposition}

Here we take the freedom to use the abbreviation $q_\cT(x)$ instead of the more
correct notation $q_\cT(\cD_\cT(x))$. The statement is proved by
straightforward computation. This is, of course, well known for $\cT=\I$ and
$\V$, where $q_\I$ is the classical cross-ratio, and $q_\V$ coincides, up to
1/6, with the Schwarzian derivative.   

The presented invariants are complete in the following sense: two data of the
same type $\cT$ are M\"obius equivalent if and only if the values of the 
corresponding invariants coincide. In particular, for any data $\cD_\cT(x)$
with the value of the invariant $q_\cT(x)=\a$ there exists a unique M\"obius
transformation bringing it into the {\it canonical representative}
$\cD_\cT=\cC_\cT(\a)$, shown in the last column of the Table \ref{tab:inv}.

Finally, we notice that the above hierarchy of M\"obius invariants may be
obtained by the following series of natural degeneration processes:\medskip\\
$\I\to\II$: let $x_4=x_3+\ep\dot x_3+O(\ep^2)$, then
\[
  q_\I(x_1,x_2,x_3,x_4)=\ep q_\II(x_1,x_2,x_3,\dot x_3);
\]
$\II\to\III$: let $x_2=x_1+\ep\dot x_1+O(\ep^2)$, then
\[
  q_\II(x_1,x_2,x_3,\dot x_3)= \ep q_\III(x_1,\dot x_1,x_3,\dot x_3);
\]
$\II\to\IV$: let $x_3=x_2+\ep\dot x_2+\ep^2\ddot x_2/2+O(\ep^3)$, then
\[
 q_\II(x_1,x_2,x_3,\dot x_3)=
   -\ep^{-1}(1+\ep q_\IV(x_1,x_2,\dot x_2,\ddot x_2));
\]
$\IV\to\V$: let 
$x_2=x_1+\ep\dot x_1+\ep^2\ddot x_1/2+\ep^3\dddot x_1/6+O(\ep^4)$, then
\[
  q_\IV(x_1,x_2,\dot x_2,\ddot x_2)=
    -\ep^{-1}(1-\ep q_\V(x_1,\dot x_1,\ddot x_1,\dddot x_1)).
\]

\section{Characterization of quadrirational maps by singularities
patterns}\label{s:class}

The following statement is an easy corollary of Theorem \ref{Th sing}.

\begin{proposition}
For a quadrirational map $F$ of type $\cT$ the polynomials
$\hat{R}(u)$, $\hat{r}(v)$ also belong to the type $\cT$. Their
roots $u_i$, $v_i$ ($i=1,\ldots,m$) can be ordered so that the
companion map $\F:(u,y)\mapsto(x,v)$ has $m$ singularities
$(u_i,y_i)$, the companion map $\F^{-1}:(x,v)\mapsto(u,y)$ has $m$
singularities $(x_i,v_i)$, and the inverse map
$F^{-1}:(u,v)\mapsto(x,y)$ has $m$ singularities $(u_i,v_i)$, all
three maps belonging to the type $\cT$, as well.
\end{proposition}
\begin{proof}
The fact that the $y$--components of the singularities $(u_i,y_i)$ coincide
with the $y$--components of the singularities $(x_i,y_i)$, follows from the
fact that they also are the roots of the polynomial $r(y)$. Indeed, the
formulas (\ref{conj1}) for the companion map yield:
\begin{equation}\label{system1 conj}
 d(y_i)u_i-b(y_i)=0,\quad -c(y_i)u_i+a(y_i)=0,\quad i=1,\ldots,m.
\end{equation}
Similarly, the $x$--components of the singularities $(x_i,v_i)$ coincide with
the $x$--components of the singularities $(x_i,y_i)$. For the proof of the
claim about the singularities of $F^{-1}$ make four 90--degree rotations. 
\end{proof}

Now we are in a position to formulate the main result which gives a complete
classification of quadrirational maps of the subclass [2:2].

\begin{theorem}\label{Th main}
Attach to each edge $x,y,u,v$ of the quadrilateral on Fig.~\ref{fig:maps} the
data of one and the same type $\cT$. The necessary and sufficient condition for
the existence of a quadrirational map $F:(x,y)\mapsto(u,v)$ of type $\cT$ with
the singularities data $(\cD_\cT(x),\cD_\cT(y))$ and such that the companion
maps $\F:(u,y)\mapsto(x,v)$ and $\F^{-1}:(x,v)\mapsto (u,y)$ have the
singularities data $(\cD_\cT(u),\cD_\cT(y))$ and $(\cD_\cT(x),\cD_\cT(v))$,
respectively, reads:
\begin{equation}\label{q's}
 q_\cT(x)=q_\cT(u)=\a,\quad q_\cT(y)=q_\cT(v)=\b,\quad \a\ne\b.
\end{equation}
If this condition is satisfied, then the corresponding map is
unique.
\end{theorem}
\begin{proof} 
Necessity follows from Theorem \ref{Th quadratic}. To prove the existence and
uniqueness, one can compute the required map as follows. The coefficients of
the polynomials $a(y),\ldots,d(y)$ have to satisfy the system of equations
consisting of the applicable ones among (\ref{system1}), (\ref{system3 1}),
(\ref{system4 1}) and (\ref{system5 1}), with $(\xi,\eta)=(x_i,y_i)$,
$\l=\l_i$, $\theta=\theta_i$ and $\k=\k_i$ (cf. (\ref{params})), and also the
analogous equations obtained by replacing $(a(y),b(y),c(y),d(y))$ by
$(d(y),-b(y),-c(y),a(y))$, with $(\xi,\eta)=(u_i,y_i)$ and the corresponding
parameters $\l=\mu_i$, $\theta=\vartheta_i$ and $\k=\kappa_i$ of the companion
map $\F$. This gives in all cases 16 homogeneous linear equations for 12
coefficients of the polynomials $a(y),\ldots,d(y)$. It is enough to study this
system for the special case when canonical representatives of the singularities
data $\cD_\cT$ are chosen on each edge $x,y,u,v$. The case by case study shows
that this system always admits a unique (up to a common factor) solution; we
give below more details for two extreme cases $\cT={\rm I}$ and $\cT={\rm V}$.
Analogously, the coefficients of the polynomials $A(x),\ldots,D(x)$ are
uniquely defined by the singularities data. The resulting {\it canonical
quadrirational maps of type $\cT$} will be listed below. 
\end{proof}

\paragraph{Remarks.} 
1) It is important to notice that the parameters $\a,\b$ attached
to the edges received in our construction an intrinsic
interpretation in terms of singularities of the map.

2) In case $\a=\b$ the above mentioned systems
also admit a unique solution, but the resulting map degenerates:
$u=y,\;v=x$, thus having no singularities at all.
\medskip

The {\it canonical map of type} $\cT$ is characterized by the singularities 
data
\begin{equation}\label{sing}
 \cD_\cT(x)=\cD_\cT(u)=\cC_\cT(\a),\quad
 \cD_\cT(y)=\cD_\cT(v)=\cC_\cT(\b) 
\end{equation}
accordingly to the Table~\ref{tab:inv}. 

\paragraph{Type I.} 
The system for the coefficients of the polynomials $a(y),\ldots,d(y)$ consists
of 8 equations (\ref{system1}) with $(\xi,\eta)=(x_i,y_i)$, $i=1,\ldots,4$, and
of 8 equations (\ref{system1 conj}). For the canonical choice (\ref{sing}) of
the singularities data it reads:
\[
\begin{array}{r}
  a_2=c_2=0,\quad b_0=d_0=0,\\
  a_0+a_1+b_1+b_2=0,\\
  c_0+c_1+d_1+d_2=0,\\
  (a_0+a_1\b)\a+(b_1\b+b_2\b^2)=0,\\
  (c_0+c_1\b)\a+(d_1\b+d_2\b^2)=0,
\end{array}\quad
\begin{array}{r}
 a_0=b_0=0,\quad c_2=d_2=0,\\
  d_0+d_1-b_1-b_2=0,\\
 -c_0-c_1+a_1+a_2=0,\\
  (d_0+d_1\b)\a-(b_1\b+b_2\b^2)=0,\\
 -(c_0+c_1\b)\a+(a_1\b+a_2\b^2)=0.
\end{array}
\]
This system admits a unique (up to a common factor) solution resulting exactly
in the map (\ref{F1}). 

\paragraph{Type II.} The canonical map is
\[
  u=\frac{\a y(x-y)}{\b x(1-y)-\a y(1-x)}\,,\quad
  v=\frac{\b x(x-y)}{\b x(1-y)-\a y(1-x)}\,.
\]
It can be obtained from (\ref{F1}) by the rescaling $\a\mapsto\ep\a$,
$\b\mapsto\ep\b$, and then sending $\ep\to0$. Upon inverting all variables
$x,y,u,v$, this map takes the form (\ref{F2}) with two simple singularities
at $(0,0)$, $(1,1)$ and one double singularity at $(\infty,\infty)$.

\paragraph{Type III.} The canonical map is 
\[
 u=\frac{\a y(x-y)}{\b x-\a y+(\b-\a)(y^2-2y)x}\,,\quad
 v=\frac{\b x(x-y)}{\b x-\a y+(\b-\a)(x^2-2x)y}\,.
\]
Upon the M\"obius transformation of variables $x=1/(x'-1)$, and similarly for
$y,u,v$, this map takes the form (\ref{F3}) with two double singularities at
$(0,0)$ and $(\infty,\infty)$. Obviously, this can be obtained from (\ref{F2})
by multiplying all variables $x,y,u,v$ by $T$ and then sending $T\to\infty$.

\paragraph{Type IV.} The canonical map is 
\[
 u=\frac{y(x-y)}{x-y+(\a-\b)xy}\,,\quad
 v=\frac{x(x-y)}{x-y+(\a-\b)xy}\,.
\]
Upon inverting all variables $x,y,u,v$ it takes the form (\ref{F4}) with a
simple singularity at $(0,0)$ and a triple singularity at $(\infty,\infty)$: 

\paragraph{Type V.} 
The system for the coefficients of the polynomials
$a(y),\ldots,d(y)$ consists of (\ref{system1}), (\ref{system3 1}),
(\ref{system4 1}) and (\ref{system5 1}), with
$(\xi,\eta)=(x_1,y_1)$, $\l=1$, $\theta=0$ and
$\k=4(\a-\b)$, and also the analogous equations obtained by
replacing $(a(y),b(y),c(y),d(y))$ by $(d(y),-b(y),-c(y),a(y))$,
with $(\xi,\eta)=(u_i,y_i)$ (with the same values of parameters
$\l$, $\theta$ and $\k$). This system reads:
\[
\begin{array}{l} 
 b_0=d_0=0,\\
 a_0+b_1=0,\\ 
 c_0+d_1=0,\\ 
 b_2+a_1=0,\\ 
 d_2+c_1=0,\\
 a_2+(\a-\b)a_0=0,\\ 
 c_2+(\a-\b)c_0=0,
\end{array}\quad
\begin{array}{l} 
 a_0=b_0=0,\\
 d_0-b_1=0,\\ 
 a_1-c_0=0,\\ 
 d_1-b_2=0,\\ 
 a_2-c_1=0,\\
 d_2+(\a-\b)d_0=0,\\ 
 c_2+(\a-\b)c_0=0,
\end{array}
\]
and admits a unique (up to a common factor) solution resulting in the map
\[
 u=\frac{y(x-y)}{x-y+(\b-\a)xy^2}\,,\quad
 v=\frac{x(x-y)}{x-y+(\b-\a)x^2y}\,.
\]
Upon inverting all variables $x,y,u,v$, it takes the form (\ref{F5}) with
quadruple singularity at $(\infty,\infty)$.

\section{3D consistency and Lax representations}\label{s:3d}

\begin{theorem}\label{Th 3D}
Let each edge of the three-dimensional cube on Fig.~\ref{fig:3D} carry the
data of the same type $\cT$, such that the M\"obius invariants $q_\cT$ of the
data on the opposite edges coincide. Attach quadrirational maps of type $\cT$
to each face of the cube, according to Theorem \ref{Th main}. Then this system
of maps is three-dimensionally consistent.
\end{theorem}
\begin{proof} 
Performing a suitable M\"obius transformation on each edge, we can assume that
all edges carry the canonical representatives of their type. Then all maps on
all faces are canonical. Consistency in this situation  follows from Theorem
\ref{th:greeks} and the fact that all canonical maps $F_\cT$ come from the
geometric construction of Sect.~\ref{s:geo}. One might prefer a more
direct proof. It can be (and has been) obtained by straightforward computations 
with the help of a symbolic computations system (we used {\it Mathematica}).
\end{proof}

It is anticipated that the 3D consistency might be also proved with the
help of a more deep study of the algebraic geometry of the maps $F_\cT$.
In this respect, an information on the blow-up and blow-down structure of 
these maps may be useful. Such information is provided in Appendices
\ref{a:blowdown}, \ref{a:blowup}. 

As postulated in \cite{ABS}, the 3D consistency of discrete 2D systems
can be taken as a definition of their integrability. In particular, one
can derive from the 3D consistency such more traditional attributes of
integrability as {\it Lax representations}. The correspondent procedure
was elaborated in \cite{ABS} for discrete models with fields on vertices,
and in \cite{SV} for discrete models with fields on edges, i.e. for
Yang--Baxter maps. The latter procedure looks as follows. Consider the
parameter-dependent Yang--Baxter map $R(\a,\b):(x,y)\mapsto(u,v)$. Suppose
that on the set $\cX$ there is an effective action of the linear group
$G=GL_N$, and that the map $R(\a,\b)$ has the following special form:
\[
 u=M(y,\b,\a)[x],\quad v=L(x,\a,\b)[y],
\]
where $L,M:\cX\times\Complex\times\Complex\mapsto GL_N$ are some matrix--valued
functions on $\cX$ depending on parameters $\a$ and $\b$, and $M[x]$ denotes
the action of $M\in G$ on $x\in\cX$. Then the following relations hold:
\begin{equation}\label{Lax}
\begin{aligned}
 L(x,\a,\l)L(y,\b,\l) &= L(v,\b,\l)L(u,\a,\l), \\
 M(y,\b,\l)M(x,\a,\l) &= M(u,\a,\l)M(v,\b,\l). 
\end{aligned}
\end{equation}
Each one of these relations serves as a Lax representation for the map
$R(\a,\b)$.

In the present case of bi-M\"obius maps (\ref{map}) the above conditions 
are obviously satisfied with $GL_2$ matrices acting on $\CP^1$ by M\"obius
transformations,
\[
 M(y)=\begin{pmatrix} a(y) & b(y) \\ c(y) & d(y) \end{pmatrix},\quad
 L(x)=\begin{pmatrix} A(x) & B(x) \\ C(x) & D(x) \end{pmatrix}.
\]
The action is in this case projective, i.e. scalar matrices (and only they) are
acting trivially, so that the resulting Lax representations are {\it \`a
priori} only projective ones, i.e. (\ref{Lax}) hold only up to a
multiplication by scalars. But in each case it is not difficult to turn  them
into proper Lax representations by finding a suitable normalization.

\paragraph{Example.} For the map (\ref{F5}) the Lax matrices obtained by the
above procedure are:
\[
 L(x,\a,\l)=\begin{pmatrix} x & \l-\a-x^2\\ 1 & -x \end{pmatrix}.
\]

\section{Concluding remarks}\label{s:n}

The notions and constructions of this paper can be generalized in various
directions. One of them is the $(d+1)$-dimensional consistency of
discrete $d$-dimensional systems, discussed in some detail in \cite{ABS}.
Another one deals with Yang-Baxter maps on more general varieties than
$\CP^1\times\CP^1$, and we briefly discuss it now.

The construction in Sect.~\ref{s:geo} leads to quadrirational maps
on $\CP^1\times\CP^1$ after a rational parametrization of linear pencils of
conics. It is not difficult to represent these maps also directly in 
terms of (non-homogeneous) coordinates on $\CP^2$. Consider a linear pencil
of conics in $\Complex^2$ or $\Real^2$ in the form
\begin{equation}\label{Cl}
  Q_\l:\quad Q(X,\l)=\<X,(\l S+T)X\>+\<\l s+t,X\>+\l\sigma+\tau=0,
\end{equation}
where $S,T$ are symmetric matrices, $s,t$ are vectors and $\sigma,\tau$ scalars.
Let $X\in Q_\a$ and $Y\in Q_\b$.
An easy computation shows that two other points of intersection of the line
$\overline{XY}$ with these conics are given by the formula
\begin{equation}\label{XY}
\begin{aligned}
  X_2&= Y+\frac{(\a-\b)(\<Y,SY\>+\<s,Y\>+\sigma)}
               {\<X-Y,(\a S+T)(X-Y)\>}(X-Y),\\
  Y_1&= X+\frac{(\a-\b)(\<X,SX\>+\<s,X\>+\sigma)}
               {\<X-Y,(\b S+T)(X-Y)\>}(X-Y).
\end{aligned}	       
\end{equation}
Of course, this formula can be rewritten in several ways using the
relations $Q(X,\a)=Q(Y,\b)=0$. However, the presented version has a remarkable
property: it does not contain $t$ and $\tau$. This allows one to interpret the
initial values $X,Y$ as {\em free} points of the plane, rather than points on
the conics $Q_\a$ and $Q_\b$, and the map $(X,Y)\to (X_2,Y_1)$ as a map
$\Complex^2\times\Complex^2\to\Complex^2\times\Complex^2$ rather than 
a map $Q_\a\times Q_\b\to Q_\a\times Q_\b$.

\begin{proposition}\label{th:free}
Formulas (\ref{XY}) define the 3D consistent map 
$\Complex^n\times\Complex^n\to\Complex^n\times\Complex^n$ for an arbitrary 
dimension $n$.
\end{proposition}
\begin{proof}
It is easy to see that all points $X_2,\dots,Z_{12}$ taking part in the
definition of 3D consistency will lie in the two-dimensional plane through
the initial points $X,Y,Z$. Therefore, the proof can be effectively reduced 
to the case $n=2$. In this case, determine the vector $t$ and the scalar $\tau$
from the equations $Q(X,\a)=Q(Y,\b)=Q(Z,\g)=0$. Obviously, these equations form
a linear system for $t,\tau$, with the determinant equal to 
zero if and only if $X,Y$ and $Z$ lie on a straight line, which is
a variety of codimension 1. For $t,\tau$ thus found consider the linear
pencil of conics (\ref{Cl}) and apply Theorem \ref{th:greeks}.
\end{proof}

Note that the map (\ref{XY}) makes sense also in the scalar case $n=1$,  where
it coincides with one of the maps from the list (\ref{F2})--(\ref{F5}).  For
example, if $S\ne 0$ then we recover the map (\ref{F2}),  so that the formula
(\ref{XY}) can be considered as its multifield generalization.

Of course, there exist more 3D consistent maps in the multi-dimensional
case. As another example we present the map
\[
  X_2= -\frac{\b}{\a}Y+\frac{\a-\b}{\a}(X^{-1}-Y^{-1})^{-1},\quad
  Y_1= -\frac{\a}{\b}X+\frac{\a-\b}{\b}(X^{-1}-Y^{-1})^{-1},
\]
which generalizes the scalar map (\ref{F3}). Here $X,Y$ may be complex matrices
or vectors, and $X^{-1}$ denotes the matrix inverse or the ``vector inverse'' 
$X^{-1}=-X/|X|^2$, respectively.

There exist also maps without a direct scalar analog, for example those
introduced in \cite{GV,NY}. It should be stressed once more that in all known
cases 3D consistent maps are quadrirational.

\section*{Appendices}\appendix
\section{Blow-down}\label{a:blowdown}

As it happens with birational maps of two-dimensional varieties, our
quadrirational maps blow up singular points to curves, and blow down
exceptional curves to points.

\begin{theorem}\label{Th blowdown}
Consider a quadrirational map $F$ of type $\cT$. Then there are exactly $m$
exceptional curves $\cC_i$, $i=1,\ldots,m$, on $\CP^1\times\CP^1$, which are
blown down by the map $F$ to the points $(u_i,v_i)$. The curves
$\cC_i=F^{-1}(u_i,v_i)$ have bidegree (1,1) and are characterized by the
following conditions.
\begin{itemize}\setlength{\itemsep}{0pt}
\item[I:] In this case $m=4$. The curve $\cC_i$ ($1\le i\le 4$) passes through
three points $(x_j,y_j)$ with $j\neq i$.
\item[II:] In this case $m=3$. For $i=1,2$, the curve $\cC_i$ passes through
two points $(x_j,y_j)$, $\,j\neq i$, with the tangent at $(x_3,y_3)$ parallel
to $(\dot x_3,\dot y_3)$. The curve $\cC_3$ passes through all three points
$(x_j,y_j)$.
\item[III:] In this case $m=2$. The curve $\cC_i$ ($i=1,2$) passes through the
both points $(x_1,y_1)$ and $(x_2,y_2)$, with the tangent at $(x_j,y_j)$,
$\,j\neq i$, parallel to $(\dot x_j,\dot y_j)$.
\item[IV:] In this case $m=2$. The curve $\cC_1$ passes through the point
$(x_2,y_2)$, and has there the 2--germ defined by $(\dot x_2,\dot y_2,\ddot
x_2,\ddot y_2)$. The curve $\cC_2$ passes through the both points $(x_j,y_j)$
with the tangent at $(x_2,y_2)$ parallel to $(\dot x_2,\dot y_2)$.
\item[V:] In this case $m=1$. The curve $\cC_1$ passes through the point
$(x_1,y_1)$ and has there the 2--germ defined by $(\dot x_1,\dot y_1,\ddot
x_1,\ddot y_1)$.
\end{itemize}
\end{theorem}
\begin{proof}
It is clear (and well known) that only singular points $(u_i,v_i)$
can be images of exceptional curves. Fix $1\le i\le m$; the
pre-image $\cC_i=F^{-1}(u_i,v_i)$ is described by the equations
\[
 \varphi(x,y,u_i)=0,\quad \Phi(y,x,v_i)=0.
\]
The function $\varphi(x,y,u_i)$ is a polynomial of degree 2 in $y$
and degree $1$ in $x$, the polynomial $\Phi(y,x,v_i)$ has degree 2
in $x$ and degree $1$ in $y$. We have:
\begin{align*}
 \varphi(x,y_i,u_i) &= (c(y_i)u_i-a(y_i))x+(d(y_i)u_i-b(y_i)) \equiv0,\\
 \Phi(y,x_i,v_i)    &= (C(x_i)v_i-A(x_i))y+(D(x_i)v_i-B(x_i)) \equiv0,
\end{align*}
therefore $\varphi(x,y,u_i)$ is divisible by $y-y_i$, and
$\Phi(y,x,v_i)$ is divisible by $x-x_i$. Set
\[
 \varphi(x,y,u_i)=(y-y_i)\g_i(x,y),\quad \Phi(y,x,v_i)=(x-x_i)\Gamma_i(x,y).
\]
Obviously, both polynomials $\g_i(x,y)$ and $\Gamma_i(x,y)$
are linear in $x$ and $y$, i.e. equating them to zero defines
curves of bidegree (1,1). We claim that the equations
\[
 \g_i(x,y)=0 \quad\text{and}\quad \Gamma_i(x,y)=0
\]
are equivalent, so that the curve $\cC_i$ is given by either of them. This is
done by demonstrating that either of these equations defines the curve with the
properties listed in the theorem. Since these properties determine the curve
uniquely, the claim will follow. This demonstration is performed for each type
$\cT$ separately. Since the argument are to a large extent analogous, we
restrict ourselves to details for two extreme cases $\cT=\I$ and $\cT=\V$
only.

For $\cT=\I$, we have: $\varphi(x_j,y_j,u)=\Phi(y_j,x_j,v)\equiv0$ for all
$1\le j\le4$. It follows that $\g_i(x_j,y_j)=\Gamma_i(x_j,y_j)=0$ for all $j\ne
i$.

For $\cT=\V$, we have to demonstrate that the following holds:
\begin{align*}
 \text{a)}\quad& \g_1(x_1,y_1)=0,\\
 \text{b)}\quad& \g_{1x}(x_1,y_1)\dot x_1+\g_{1y}(x_1,y_1)\dot y_1=0,\\
 \text{c)}\quad& 2\g_{1xy}(x_1,y_1)\dot x_1\dot y_1+
  \g_{1x}(x_1,y_1)\ddot x_1+\g_{1y}(x_1,y_1)\ddot y_1=0,
\end{align*}
and the same with the replacement of $\g_1$ by $\Gamma_1$. We do this for
$\g_1$ only, since for $\Gamma_1$ everything is similar. Claim a) is proved by
the following computation:
\begin{align*}
 \g_1(x_1,y_1) 
  &= \lim_{y\to y_1}\frac{\varphi(x_1,y,u_1)}{y-y_1}=\varphi_y(x_1,y_1,u_1)\\
  &= (c'(y_1)x_1+d'(y_1))u_1-(a'(y_1)x_1+b'(y_1))\\
  &= -\l_1(c(y_1)u_1-a(y_1))=0.
\end{align*}
To check claim b), we proceed as follows:
\begin{align*}
 \g_{1x}(x_1,y_1) 
 &= \lim_{y\to y_1}\frac{\varphi_x(x_1,y,u_1)}{y-y_1}
  =\varphi_{xy}(x_1,y_1,u_1)\\
 &= c'(y_1)u_1-a'(y_1),\\
 \g_{1y}(x_1,y_1) 
 &= \lim_{y\to y_1}\frac{(y-y_1)\varphi_y(x_1,y,u_1)-\varphi(x_1,y,u_1)}
      {(y-y_1)^2}\\
 &= \tfrac12\varphi_{yy}(x_1,y_1,u_1)\\
 &= \tfrac12(c''(y_1)x_1+d''(y_1))u_1-(1/2)(a''(y_1)x_1+b''(y_1))\\
 & = -\l_1(c'(y_1)u_1-a'(y_1))
     -\tfrac12\l_1^{3/2}\theta_1(c(y_1)u_1-a(y_1))\\
 & = -\l_1(c'(y_1)u_1-a'(y_1)).
\end{align*}
To prove the claim c), we derive further:
\begin{align*}
 \g_{1xy}(x_1,y_1) 
 &= \lim_{y\to y_1}\frac{(y-y_1)\varphi_{xy}(x_1,y,u_1)-\varphi_x(x_1,y,u_1)}
      {(y-y_1)^2}\\
 &= \tfrac12\varphi_{xyy}(x_1,y_1,u_1)\\
 &= \tfrac12(c''(y_1)u_1-a''(y_1))\\
 &= -\tfrac12\l_1^{1/2}\theta_1(c'(y_1)u_1-a'(y_1))
    -\tfrac12\k_1\l_1^{-1}(c(y_1)u_1-a(y_1))\\
 &= -\tfrac12\l_1^{1/2}\theta_1(c'(y_1)u_1-a'(y_1)).
\end{align*}
Now the claim c) becomes equivalent to:
\[
 \l_1^{1/2}\theta_1\dot x_1\dot y_1=\ddot x_1-\l_1\ddot y_1\,,
\]
which, in turn, is equivalent to the second formula in (\ref{params}). 
\end{proof}

\section{Blow-up}\label{a:blowup}

Consider a quadrirational map $F$ of type $\cT=\I,\II,\III,\IV$
or $\V$. Then the point $(x_i,y_i)$ ($i=1,\ldots,m$) is blown up to
the curve $\cD_i$ on $\CP^1\times\CP^1$ defined as the
curve of bidegree (1,1) satisfying the conditions analogous to
those listed in Theorem \ref{Th blowdown}. We want to give a more
complete analytic picture of the blow--up process.

\begin{theorem}\label{Th blowup}
Let $(x_i,y_i)$ be a singular point of the map $F$, such that $x_i,y_i$ are
roots of $R(x),r(y)$, respectively, of multiplicity $1\le\ell\le4$. Let
$\cG^{(\ell-1)}_i=(\dot x_i,\dot y_i,\ldots,x^{(\ell-1)}_i,y^{(\ell-1)}_i)$ 
be the corresponding singularity data of the map $F$ at the point $(x_i,y_i)$. 
\begin{itemize}\setlength{\itemsep}{0pt}
\item If $(x(t),y(t))\to(x_i,y_i)$ along a curve with a definite 
$(\ell-1)$--germ different from $\cG^{(\ell-1)}_i$, then 
$F(x(t),y(t))\to(u_i,v_i)$.
\item If $(x(t),y(t))\to(x_i,y_i)$ along a curve with a definite $\ell$--germ 
whose $(\ell-1)$--germ coincides with $\cG^{(\ell-1)}_i$, then the images 
$F(x(t),y(t))$ have a definite limiting point lying on $\cD_i$. Different 
$\ell$--germs correspond to different limiting points on $\cD_i$. 
\end{itemize}
\end{theorem}
\begin{proof} 
Note that $\ell$--germs with a fixed $(\ell-1)$--germ $\cG^{(\ell-1)}_i$ are
characterized by one number from $\CP^1$ (e.g., for $\ell=1$ --- by the slope
$\dot{x}(0)/\dot{y}(0)$, for $\ell=2$ --- by the curvature
$(\ddot{x}(0)\dot{y}(0)-\dot{x}(0)\ddot{y}(0))/(\dot{x}(0)\dot{y}(0))^{3/2}$,
etc.) Therefore the second statement of the theorem shows the isomorphism of
$\cD_i$ with $\CP^1$.  We will give the details of computations for the proof
of the second statement in the simplest case $\ell=1$ only. On a line
$\{(x,y)=(x_i+\mu t,y_i+\nu t)\}$ through $(x_i,y_i)$ the values of
$(u,v)=F(x,y)$ have a well--defined limit as $t\to 0$:
\[
 u\to\frac{a(y_i)\mu+(a'(y_i)x_i+b'(y_i))\nu}
          {c(y_i)\mu+(c'(y_i)x_i+d'(y_i))\nu},\quad
 v\to\frac{A(x_i)\nu+(A'(x_i)y_i+B'(x_i))\mu}
          {C(x_i)\nu+(C'(x_i)y_i+D'(x_i))\mu}.
\]
We have for these limit points:
\begin{align*}
 \frac{\mu}{\nu} 
  &= \frac{(c'(y_i)x_i+d'(y_i))u-(a'(y_i)x_i+b'(y_i))}{-c(y_i)u+a(y_i)}\\ 
  &= \frac{-C(x_i)v+A(x_i)}{(C'(x_i)y_i+D'(x_i))v-(A'(x_i)y_i+B(x_i))}.
\end{align*}
So, it remains to prove the following claim:
\begin{itemize}\setlength{\itemsep}{0pt}
\item If $x_i,y_i$ are simple roots of $R(x),r(y)$, respectively,
then the equation of the curve $\cD_i$ reads:
\begin{multline}\label{uv curve}
 (c(y_i)u-a(y_i))(C(x_i)v-A(x_i)) \\
   = \Big((c'(y_i)x_i+d'(y_i))u-(a'(y_i)x_i+b'(y_i))\Big) \\
     \times\Big((C'(x_i)y_i+D'(x_i))v-(A'(x_i)y_i+B(x_i))\Big).
\end{multline}
\end{itemize}
According to the proof of Theorem \ref{Th blowdown}, applied with
$F^{-1}$ instead of $F$, the equations of $\cD_i$ are given by
\[
 \hat{\varphi}(u,v,x_i)/(v-v_i)=0 \quad\text{or}\quad
 \hat\Phi(v,u,y_i)/(u-u_i)=0.
\]
To show that these equations are equivalent to (\ref{uv curve}), 
recall that, for instance,
\begin{align*}
 \hat\Phi(v,u,y)
  &= (\hat{C}(u)y+\hat{D}(u))v-(\hat{A}(u)y+\hat{B}(u))\\
  &= \frac{1}{r(y)}(q(y,u)v-p(y,v))\\
  &= \frac{(c(y)u-a(y))^2}{r(y)}\Big((C(X)y+D(X))v-(A(X)y+B(X))\Big),
\end{align*}
where
\[
 X=X(y,u)=\frac{d(y)u-b(y)}{-c(y)u+a(y)}.
\]
Recall also that for $y=y_i$ we have:
\[
 X(y_i,u)=\frac{d(y_i)u-b(y_i)}{-c(y_i)u+a(y_i)}\equiv x_i\,.
\]
Therefore:
\begin{multline*}
 \hat\Phi(v,u,y_i)=\frac{(c(y_i)u-a(y_i))^2}{r'(y_i)}\,(C(x_i)v-A(x_i))\\
  +\frac{(c(y_i)u-a(y_i))^2}{r'(y_i)}
  \Big((C'(x_i)y_i+D'(x_i))v-(A'(x_i)y_i+B(x_i))\Big)X_y(y_i,u).
\end{multline*}
It is easy to calculate that:
\[
 X_y(y_i,u)=\frac{(c'(y_i)x_i+d'(y_i))u-(a'(y_i)x_i+b'(y_i))}
   {-c(y_i)u+a(y_i)}\,.
\]
Noting also that $c(y_i)u-a(y_i)=c(y_i)(u-u_i)$, we find finally
the equations of the curve $\cD_i$ in the required form (\ref{uv
curve}):
\begin{multline*}
 \frac{r'(y_i)}{c(y_i)}\cdot\frac{\hat\Phi(v,u,y_i)}{u-u_i}
  =(c(y_i)u-a(y_i))(C(x_i)v-A(x_i))\\
  -\Big((c'(y_i)x_i+d'(y_i))u-(a'(y_i)x_i+b'(y_i))\Big)\\
  \times\Big((C'(x_i)y_i+D'(x_i))v-(A'(x_i)y_i+B(x_i))\Big).
\end{multline*}
This proves the second statement of the theorem for $\ell=1$. Calculation 
for $\ell>1$ are similar, as well as the proof of the first statement. 
\end{proof}

\section{Subclasses [1:2] and [1:1]}\label{a:12}

In this Appendix we discuss briefly maps of the subclasses [1:2] and
[1:1] (cf. Sect.~\ref{s:BD}). According to Lemma \ref{lemma semilinear}, 
these maps are described by the polynomial systems of the form
\begin{equation}\label{BD1}
 \Phi(y,x,v)=f(y,v)x+g(y,v)=0,\quad 
 \hat\Phi(y,u,v)=\hat f(y,v)u+\hat g(y,v)=0, 
\end{equation} 
where $f(y,v)=f_1yv+f_2y+f_3v+f_4$ and analogously for $g,\hat f,\hat g$. In
order to find a canonical form of this system with respect to the
transformations from $(\cM\ddot{o}b)^4$, notice that M\"obius transformations 
$ x\to(\k x+\mu)/(\l x+\nu)$ lead to $f\to\k f+\l g$, 
$g\to\mu f+\nu g$.
They can be used to turn the polynomials $f,g$ into reducible ones.
It is easy to see that this is possible if and only if the characteristic
polynomial
\[
  \det\left(\k\begin{pmatrix} f_1 & f_2\\ f_3 & f_4 \end{pmatrix}
           +\l\begin{pmatrix} g_1 & g_2\\ g_3 & g_4 \end{pmatrix}\right)
\] 
has two distinct eigenvalues $(\k:\l),(\mu:\nu)\in\CP^1$. If this is true
also for the second equation of the system (\ref{BD1}) (i.e. if the polynomials 
$\hat f,\hat g$ can be made reducible by a M\"obius transformation in $u$), 
then the map (\ref{BD1}) can be brought into the form
\begin{equation}\label{BDred}
\begin{aligned}
  x(y-y_0)(v-v_0)+(y-y_1)(v-v_1)&=0,\\
  u(y-y_2)(v-v_2)+(y-y_3)(v-v_3)&=0, 
\end{aligned} 
\end{equation}
with some $y_i,v_i\in\CP^1$. If each quadruple $(y_i)$ and  $(v_i)$ contains 
four distinct points then suitable M\"obius transformations in $y$ and $v$ 
bring these quadruples into $(\infty,0,\b,1)$ and $(\a,1,\infty,0)$,
respectively, and the system (\ref{BDred}) -- into the form 
\[
  x(v-\a)-y(v-1)=0,\quad u(y-\b)-v(y-1)=0.
\]
This results in the map
\begin{equation}\label{12.1}
G_{\a,\b}:\quad  
u=\frac{(\a x-y)(y-1)}{(x-y)(y-\b)},\quad v=\frac{\a x-y}{x-y}.
\end{equation}
This is the canonical form of the generic map from the subclass [1:2]. 
Several degeneracies can happen, due to coincidence of some points $(y_i)$
or $(u_i)$ in (\ref{BDred}), as well as the case of double eigenvalue
$(\k:\l)=(\mu:\nu)$, for one or both equations in (\ref{BD1}).
We leave the analysis of these degenerate cases to the interested 
reader. 

\begin{figure}[htbp]
\begin{center}\includegraphics[width=12cm]{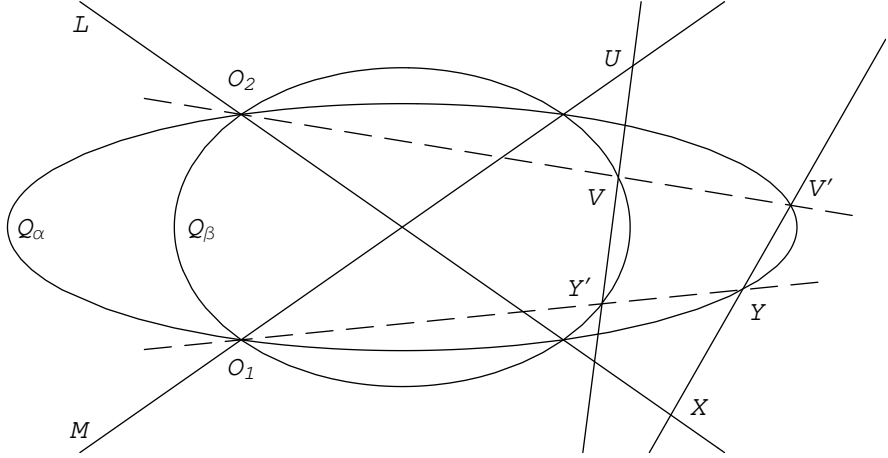}\end{center}
\caption{The geometric interpretation of the map from [1:2] subclass}
\label{fig:greeks_deg}
\end{figure}

Surprisingly, the geometric interpretation of the map (\ref{12.1}) is
somewhat more complicated than the one presented on Fig.~\ref{fig:greeks}. The
construction actually includes {\it three} conics from the linear pencil:  the
pair of nondegenerate conics $Q_\a$, $Q_\b$, and a degenerate conic which is a
pair of straight lines $L\cup M$. Fix two points $O_1,O_2$ of the base locus
of the pencil (consisting in the present generic case of four points). Choose
the degenerate conic $L\cup M$ so that neither of these lines contains both of
$O_i$. The projections from the points $O_1$, $O_2$ establish two different
isomorphisms of the nondegenerate conics of the pencil. The geometric map
described by the formulas (\ref{12.1}) is $\cF:L\times Q_\a\ni(X,Y)\mapsto
(U,V)\in M\times Q_\b$ constructed as follows:
\begin{itemize}\setlength{\itemsep}{0pt}
\item The point $V'$ is the intersection point of $Q_\a\cap\overline{XY}$ 
differerent from $Y$, and $V\in Q_\b$ is obtained from $V'\in Q_\a$  via
projection from $O_2$.
\item The point $Y'\in Q_\b$ is obtained from $Y\in Q_\a$ via  projection from
$O_1$, and $U$ is the intersection of $M$ with  $\overline{VY'}$.
\end{itemize}
This construction is illustrated on Fig.~\ref{fig:greeks_deg} (the two
identifying projections are shown by the dashed lines). The map (\ref{12.1})
results upon fixing the objects on this figure as follows. Use
non-homogeneous coordinates $(W_1,W_2)$ on the affine part $\Complex^2$ of 
$\CP^2$. The base locus of
the pencil consists of the four points  $O_1=(0,0)$, $O_2=(1,0)$, $(0,1)$ and
$(1,1)$. The equations and the parametrizations of the lines and conics 
are:
\begin{align*}
    L:\quad & W_1+W_2=1, & (W_1,W_2)&=\Bigl(\frac{1}{1+x},
                                            \frac{x}{1+x}\Bigr),\\
    M:\quad & W_1=W_2,   & (W_1,W_2)&=\Bigl(\frac{1}{1+u},
                                            \frac{1}{1+u}\Bigr),\\
 Q_\a:\quad & W_2(W_2-1)=\a W_1(W_1-1),& (W_1,W_2)&=
      \Bigl(\frac{y-\a}{y^2-\a},\frac{y(y-\a)}{y^2-\a}\Bigr)\\
  &&&=\Bigl(\frac{v-\a}{v^2-\a},\frac{\a(v-1)}{v^2-\a}\Bigr).
\end{align*}    
So, the points $X\in L$ and $Y\in Q_\a$ are parametrized by the slope $x$, 
resp. $y$, of the line $\overline{O_1X}$, resp. $\overline{O_1Y}$, while
the points $U\in M$ and $V\in Q_\b$ are parametrized by the slope $-u$, 
resp. $-v$, of the line $\overline{O_2U}$, resp. $\overline{O_2V}$.

Note that by $\a=\b$ only two conics ($Q_\b$ and the degenerate one $L\cup M$)
are involved, and the projections from $O_1$ and $O_2$  drop out from the
construction. Correspondingly, the map (\ref{12.1}) with $\a=\b$ may be
obtained from the map of Sect.~\ref{s:geo} by considering two different
parametrizations for $X,Y$ and for $U,V$, as above, and by letting one of the
two conics  degenerate. More precisely, the map (\ref{F1}) under the M\"obius
transformations $u\mapsto\a(u-1)/(u-\a)$, $v\mapsto\b(v-1)/(v-\b)$ turns into
\[
\hat{F}_{\a,\b}:\quad u=\frac{(\b x-\a y)(y-1)}{(x-y)(y-\b)},\quad
 v=\frac{(\b x-\a y)(x-1)}{(x-y)(x-\a)},
\]
and setting here $\a=1$ (i.e. letting $Q_\a$ degenerate) leads to the same map
as setting $\a=\b$ in (\ref{12.1}).

Maps of the subclass [1:2] can be included into 3D-consistent
systems, the consistency being governed by matching singularities 
along the edges of the elementary cube, like in Theorem \ref{Th 3D}. 
Not going into details, we provide
the reader just with a generic instance of such a construction. Let both
facets of the cube on Fig.~\ref{fig:3D} parallel to the coordinate plane $xy$ 
carry the [1:2] map 
\[
\hat{G}_{\a,\b}:\quad u=\frac{(\b x-\a y)(y-1)}{(x-y)(y-\b)},\quad
 v=\frac{(\b x-\a y)}{\a(x-y)},
\]
both facets parallel to the coordinate plane $xz$ carry the [1:2] 
map $\hat{G}_{\a,\g}$, while the facets parallel to the coordinate
plane $yz$ carry the [2:2] maps $\hat{F}_{\b,\g}$ and $\hat{F}_{\b/\a,\g/\a}$
(the left and the right ones, respectively). Then this system of maps
is 3D consistent. (Note that $\hat{G}_{\a,\b}$ coincides with
$G_{\b/\a,\b}$ from (\ref{12.1}) up to a scaling of the variable $u$).

Finally, note that the [1:1] maps may be interpreted in the same way, with two
degenerate conics (two pairs of lines). A generic map of this subclass is
a parameter-free one: $u=(x+y-1)/y,\; v=(x+y-1)/x$.
This map is 3D consistent, so that its companion
$u=xv/(1-x+xv), \; y=1-x+xv$ is a Yang-Baxter map (example given in
 \cite{D,O2}).


\end{document}